\documentclass{article}[12pt]
\usepackage{hyperref}
\usepackage{amsmath,amssymb,latexsym,amsthm}
\usepackage{url}
\usepackage{graphicx,color}
\usepackage{mathdots}
\usepackage[ruled,linesnumbered]{algorithm2e}
\usepackage[margin=1.2in]{geometry}
\usepackage{caption}
\usepackage{mathtools}
\usepackage{soul}
\usepackage{authblk}

\newtheorem{thm}{Theorem}
\newtheorem{lemma}[thm]{Lemma}
\newtheorem{prop}[thm]{Proposition}

\newtheorem{coro}[thm]{Corollary}
\newtheorem{remark}[thm]{Remark}

\DeclarePairedDelimiter\norm{\lVert}{\rVert}
\DeclareMathOperator*{\dif}{d\!}

\def\p{\partial}

\graphicspath{{figures/}}

\title{Linearized Boundary Control Method for Density Reconstruction in Acoustic Wave Equations}

\author[1]{Lauri Oksanen \thanks{lauri.oksanen@helsinki.fi}}
\author[2]{Tianyu Yang
\thanks{yangti27@msu.edu}}
\author[3]{Yang Yang
\thanks{yangy5@msu.edu}}
\affil[1]{Department of Mathematics and Statistics, University of Helsinki}
\affil[2]{Department of Computational Mathematics, Science and Engineering, Michigan State University}
\affil[3]{Department of Computational Mathematics, Science and Engineering, Michigan State University}

\date{}

\begin{document}

\maketitle

\abstract{We develop a linearized boundary control method for the inverse boundary value problem of determining a density in the acoustic wave equation.
The objective is to reconstruct an unknown perturbation in a known background density from the linearized Neumann-to-Dirichlet map.
A key ingredient in the derivation is a linearized Blagove\u{s}\u{c}enski\u{ı}’s identity with a free parameter.
When the linearization is at a constant background density, we derive two reconstructive algorithms with stability estimates based on the boundary control method.
When the linearization is at a non-constant background density, we establish an increasing stability estimate for the recovery of the density perturbation.
The proposed reconstruction algorithms are implemented and validated with several numerical experiments to demonstrate the feasibility.}

\section{Introduction}

The paper is concerned with the linearized inverse boundary value problem (IBVP) for the acoustic wave equation with a potential. The goal is to derive uniqueness, stability estimates and reconstruction procedures to numerically compute a small perturbation of a certain parameter in the wave equation from the knowledge of infinite boundary data. This parameter reduces to the slowness (i.e, reciprocal of the wave speed) in the absence of the potential.

\medskip
\textbf{Formulation.}
We begin with the formulation of the inverse boundary value problem for the wave equation. Let $T>0$ be a constant and $\Omega\subset\mathbb{R}^n$ be a bounded open subset with smooth boundary $\partial\Omega$. 
Consider the following boundary value problem for the acoustic wave equation with potential: 
\begin{equation} \label{eq:bvp}
\left\{
\begin{array}{rcl}
  \rho(x) \partial^2_t u(t,x) - \Delta u(t,x) + q(x) u(t,x) & = & 0, \quad\quad\quad \text{ in } (0,2T) \times \Omega \\
 \partial_\nu u  & = & f, \quad\quad\quad \text{ on } (0,2T) \times \partial\Omega \\ 
 u(0,x) = \partial_t u(0,x) & = & 0 \quad\quad\quad\quad x \in \Omega.
\end{array}
\right.
\end{equation}
Here, $\rho(x) \in C^\infty(\overline{\Omega})$ is a positive smooth function that is strictly positive, $q(x)\in C^\infty(\overline{\Omega})$ is a real-valued function referred to as the \textit{potential}. We write the wave solution as $u = u^f(t,x)$ whenever it is necessary to specify the Neumann data.

Given $f\in C^\infty_c((0,2T)\times\partial\Omega)$, the well-posedness of this problem is ensured by the standard theory for second order hyperbolic partial differential equations~\cite{evans1998partial}.
As a result, the following \textit{Neumann-to-Dirichlet map (ND map)} is well defined:
\begin{equation} \label{eq:NDmap}
\Lambda_{\rho} f := u^f|_{(0,2T)\times \partial\Omega}.
\end{equation}

Throughout the paper, we will fix a \textit{known} potential $q$ and only study the dependence of the ND map on the parameter $\rho$. This dependence is indicated by the subscript.
The inverse boundary value problem (IBVP) we are interested in concerns recovery of $\rho(x)$ from knowledge of the ND map $\Lambda_{\rho}$, that is, to invert the \textit{parameter-to-data map} $\rho\mapsto \Lambda_{\rho}$.

\medskip
\textbf{Literature.} The IBVP has been studied in the mathematical literature for a long time, and most of the results were obtained in the absence of the potential. When $q\equiv 0$, the parameter $\rho(x)$ is related to the wave speed $c(x)$ by $\rho(x)=c^{-2}(x)$. In this circumstance, the IBVP aims to recover the spatial distribution of the wave speed $c$. 
For $q\equiv 0$ and $c=c(x)$, Belishev~\cite{belishev1988approach} proved that $c$ (hence $\rho$) is uniquely determined using the boundary control (BC) method combined with Tataru's unique continuation result~\cite{tataru1995unique, tataru1999unique}.
The method has since been extended to 
many wave equations.
We mention~\cite{belishev1992reconstruction} for a generalization to Riemannian manifolds, and~\cite{MR1331288} for a result covering all symmetric time-independent lower order perturbations of the wave operator. Non-symmetric, time-dependent and matrix-valued lower order perturbations were recovered in~\cite{eskin2007inverse} \cite{MR1784415}, and~\cite{MR3880231}, respectively. For a review of the BC method, we refer to~\cite{MR2353313, MR1889089}.

As for the stability, it can be proven that the IBVP to recover the wave speed $c$ is H\"older stable under suitable geometric assumptions~\cite{stefanov1998stability, stefanov2005stable}, even when the speed is anisotropic (hence represented by a Riemannian metric). 
On the other hand, a low-pass version of $c$ can be recovered with Lipschitz-type stability~\cite{liu2016lipschitz}.

The BC method has been implemented numerically to reconstruct the wave speed $c$ in~\cite{belishev1999dynamical}, and subsequently in~\cite{belishev2016numerical, de2018recovery, pestov2010numerical, yang2021stable}. 
The implementations~\cite{belishev1999dynamical, belishev2016numerical, de2018recovery} involve solving unstable control problems, whereas~\cite{pestov2010numerical, yang2021stable} are based on solving stable control problems but with target functions exhibiting exponential growth or decay. The exponential behaviour leads to instability as well. On the other hand, the linearized approach introduced in the present paper is stable. 
It should be mentioned that the BC method can be implemented in a stable way in the one-dimensional case, see~\cite{korpela2018discrete}.
For an interesting application of a variant of the method in the one-dimensional case, see~\cite{blaasten2019blockage} on detection of blockage in networks.

\medskip
\textbf{Formal Linearization.}
We are interested in the linearized IBVP in this paper. To this end, let us formally linearize the parameter-to-data map $\rho\mapsto \Lambda_\rho$. Recall that the potential $q$ is known and fixed. We write
$$
\rho(x) = \rho_0(x) + \epsilon\dot{\rho}(x), \quad\quad u(t,x) = u_0(t,x) + \epsilon \dot{u}(t,x)
$$
where $\rho_0$ is a known background potential and $u_0$ is the background solution, $\dot{\rho}\in C_c^\infty(\Omega)$ is a compactly supported smooth perturbation.
Equating the $O(1)$-terms gives
\begin{equation} \label{eq:unperturbedBVP}
\left\{
\begin{array}{rcll}
\rho_0(x) \partial^2_t u_0(t,x) - \Delta u_0(t,x) + q(x) u_0(t,x) & = & 0, &\quad \text{ in } (0,2T) \times \Omega \\
 \partial_\nu u_0  & = & f, &\quad \text{ on } (0,2T) \times \partial\Omega \\ 
 u_0(0,x) = \partial_t u_0(0,x) & = & 0, &\quad x \in \Omega.
\end{array}
\right.
\end{equation}
Equating the $O(\epsilon)$-terms gives
\begin{equation} \label{eq:linearizedBVP}
\left\{
\begin{array}{rcll}
\rho_0(x) \partial^2_t \dot{u}(t,x) - \Delta \dot{u}(t,x) + q(x) \dot{u}(t,x) & = & -\dot{\rho}\partial^2_t u_0(t,x), &\quad \text{ in } (0,2T) \times \Omega \\
 \partial_\nu \dot{u}  & = & 0, &\quad \text{ on } (0,2T) \times \partial\Omega \\ 
 \dot{u}(0,x) = \partial_t \dot{u}(0,x) & = & 0 &\quad x \in \Omega.
\end{array}
\right.
\end{equation}

Correspondingly, we write the ND map as $\Lambda_{\rho}=\Lambda_{\rho_0} + \epsilon \dot{\Lambda}_{\dot{\rho}}$, where $\Lambda_{\rho_0}$ is the ND map for the unperturbed boundary value problem~\eqref{eq:unperturbedBVP}, and $\dot{\Lambda}_{\dot{\rho}}$ is defined as
\begin{equation} \label{eq:linearizedNDmap}
\dot{\Lambda}_{\dot{\rho}}: f \mapsto \dot{u}|_{(0,2T)\times\partial\Omega}.
\end{equation}
Note that the unperturbed problem~\eqref{eq:unperturbedBVP} can be explicitly solved to obtain $u_0$ and $\Lambda_{\rho_0}$, since $\rho_0$ and $q$ are known. As before, we will write $\dot{u} = \dot{u}^f := \dot{\Lambda}_{\dot{\rho}} f$ if it is necessary to specify the Neumann data $f$.
The linearized IBVP concerns recovery of a compactly-supported smooth perturbation $\dot{\rho}\in C^\infty_c(\Omega)$ near a known $0<\rho_0\in C^\infty(\overline{\Omega})$ from the data $\dot{\Lambda}_{\dot{\rho}}$.

\medskip
\textbf{Contribution of the Paper.} 
The major contribution of this paper consists of novel ideas to tackle linearized acoustic IBVPs as well as several results regarding the uniqueness, stability, and reconstructive algorithms. These include
\begin{itemize}
    \item A linearized Blagove\u{s}\u{c}enski\u{ı}’s identity with a free parameter. The Blagove\u{s}\u{c}enski\u{ı}’s identity plays a central role in the boundary method by bridging boundary measurement with inner products of waves. For the linearized IBVP, the authors' earlier work~\cite{oksanen2022linearized} derived a version of the linearized Blagove\u{s}\u{c}enski\u{ı}’s identity with a free parameter in the presence of a potential. The free parameter enlarges the class of testing functions that can be used to probe the unknown parameter, resulting in improved stability and reconstruction. In this paper, we derived another version of the linearized Blagove\u{s}\u{c}enski\u{ı}’s identity with a free parameter in the presence of a density, see Proposition~\ref{thm:id}. This identity forms the foundation for the derivation of other results.
    \item Multiple reconstruction formulae and algorithms to recover $\dot{\rho}$. Specifically: (1) For constant $\rho_0$ and $q\equiv 0$, we derive a reconstruction formula for $\dot{\rho}$ in Algorithm~\ref{alg:constq}. This algorithm is numerically implemented and validated in 1D with quantitative assessment of the accuracy. Moreover, A pointwise stability estimate for the Fourier transform of $\dot{\rho}$ is established in Theorem~\ref{thm:stab1}. (2) For constant $\rho_0$ and $q\not\equiv 0$, we derive another reconstruction formula for $\dot{\rho}$ in Algorithm~\ref{alg:varq}.
    \item An increasing stability estimate. For variable $\rho_0\in C^\infty(\overline{\Omega})$, we prove a stability estimate for the reconstruction of $\dot{\rho}$ in Theorem~\ref{thm:stab}. This estimate contains the free parameter in the linearized Blagove\u{s}\u{c}enski\u{ı}’s identity. The stability is a blend of a H\"older-type stability and a log-type stability. However, as the free parameter increases, the log-type decreases, leading to a nearly H\"older-type stability. This phenomenon is known as the increasing stability, and has been studied in the frequency domain for the Helmholtz equation~\cite{cheng2016increasing, hrycak2004increased, isakov2010increasing, isakov2014increasing, isakov2016increasing, isakov2018increasing, kow2021optimality, nagayasu2013increasing}. In this paper, we make use of the linearized Blagove\u{s}\u{c}enski\u{ı}’s identity with a free parameter to establish an increasing stability result for the IBVP in the time domain, see Theorem~\ref{thm:stab}. An interesting observation is that the free parameter plays the role of the frequency for the probing test functions.
\end{itemize}

\medskip
\textbf{Paper Structure.} 
The paper is organized as follows. Section~\ref{sec:prelim} reviews fundamental concepts and results in the BC method. Section~\ref{sec:id} is devoted to the proof of an integral identity that is essential to the development of our linearized BC method. Section~\ref{sec:stabandrecon} establishes several stability estimates and reconstructive algorithms for the linearized IBVP, which are the main results of the paper. Section~\ref{sec:experiment} consists of numerical implementation of a reconstruction formula as well as multiple numerical experiments for a proof-of-concept validation.

\section{Preliminaries} \label{sec:prelim}

Introduce some notations: 
Given a function $u(t,x)$, we write $u(t)=u(t,\cdot)$ for the spatial part as a function of $x$.
Introduce the time reversal operator $R: L^2([0,T]\times\partial\Omega) \rightarrow L^2([0,T]\times\partial\Omega)$,
\begin{equation} \label{eq:R}
Ru(t,\cdot):=u(T-t,\cdot) ,\quad\quad 0<t<T;
\end{equation}
and the low-pass filter $J: L^2([0,2T]\times\partial\Omega) \rightarrow L^2([0,T]\times\partial\Omega)$
\begin{equation} \label{eq:J}
Jf(t,\cdot):=\frac{1}{2}\int^{2T-t}_t f(\tau,\cdot) \,d\tau,\quad\quad 0<t<T.
\end{equation}
Let $\mathcal{T}_D$ and $\mathcal{T}_N$ be the Dirichlet and Neumann trace operators respectively, that is,
$$
\mathcal{T}_D u(t,\cdot) = u(t,\cdot)|_{\partial\Omega}, \quad\quad\quad \mathcal{T}_N u(t,\cdot) = \partial_\nu u(t,\cdot)|_{\partial\Omega}.
$$
We write $P_T: L^2([0,2T]\times\partial\Omega) \rightarrow L^2([0,T]\times\partial\Omega)$
for the orthogonal projection via restriction. Its adjoint operator $P^*_T: L^2([0,T]\times\partial\Omega) \rightarrow L^2([0,2T]\times\partial\Omega)$ is the extension by zero.

Introduce the \textit{connecting operator}
\begin{equation} \label{eq:K}
K:= J \Lambda_{\rho} P^\ast_T - R \Lambda_{\rho,T} R J P^\ast_T
\end{equation}
where $\Lambda_{\rho,T}f := P_T \Lambda_{\rho} 
 P^*_T f$.
Then the following Blagove\u{s}\u{c}enski\u{ı}’s identity holds~\cite{bingham2008iterative,blagoveshchenskii1967inverse,de2018exact,oksanen2013solving}. 

\bigskip
\begin{prop} \label{thm:waveinner}
Let $u^{P^*_T f}, u^{P^*_T h}$ be the solutions of \eqref{eq:bvp} with the Neumann traces $f, h \in L^2((0,T)\times\partial\Omega)$, respectively. Then
\begin{equation} \label{eq:ufuh}
(\rho u^{P^*_T f}(T), u^{P^*_T h}(T))_{L^2(\Omega)} = (f,Kh)_{L^2((0,T)\times\partial\Omega)} = (Kf,h)_{L^2((0,T)\times\partial\Omega)}.
\end{equation}
\end{prop}

\begin{proof}
We begin with the additional assumption that $f, h \in C^\infty_c((0,T)\times\partial\Omega)$.
Define 
$$
I(t,s) := (\rho u^{P^*_T f}(t), u^{P^*_T h}(s))_{L^2(\Omega)}.
$$
We compute
\begin{align}
 & (\partial^2_t - \partial^2_s) I(t,s) \nonumber \\
= & ((\Delta-q) u^{P^*_T f}(t), u^{P^*_T h}(s))_{L^2(\Omega)} - (u^{P^*_T f}(t),(\Delta-q) u^{P^*_T h}(s))_{L^2(\Omega)} \nonumber \\
= & (P^*_T f(t), \Lambda_{\rho} P^*_T h(s))_{L^2(\partial\Omega)} - (\Lambda_{\rho} P^*_T f(t), P^*_T h(s))_{L^2(\partial\Omega)}, \label{eq:id1}
\end{align}
where the last equality follows from the integration by parts.
On the other hand, $I(0,s)=\partial_t I(0,s) = 0$ since $u^{P^*_T f}(0,x)=\partial_t u^{P^*_T f}(0,x) = 0$. Solve the inhomogeneous $1$D wave equation \eqref{eq:id1} together with these initial conditions to obtain 
\begin{align*}
I(T,T) & = \frac{1}{2} \int^T_0 \int^{2T-t}_{t} \left[ (P^*_T f(t), \Lambda_{\rho} P^*_T h(\sigma))_{L^2(\partial\Omega)} - (\Lambda_{\rho} P^*_T f(t), P^*_T h(\sigma))_{L^2(\partial\Omega)} \right] \,d\sigma dt \vspace{1ex}\\
 & = \int^T_0  [ ( P^*_T f(t) ,  \frac{1}{2}\int^{2T-t}_{t} \Lambda_{\rho} P^*_T h(\sigma) \,d\sigma)_{L^2(\partial\Omega)} - ( \Lambda_{\rho} P^*_T f(t),  \frac{1}{2}\int^{2T-t}_{t} P^*_T h(\sigma) \,d\sigma)_{L^2(\partial\Omega)} ] \,dt \vspace{1ex} \\
 & = (f, J \Lambda_{\rho} P^*_T h)_{L^2((0,T)\times\partial\Omega)} - (P_T(\Lambda_{\rho} P^*_T f),J P^*_T h)_{L^2((0,T)\times\partial\Omega)}.
\end{align*}
Using the relations $P_T(\Lambda_{\rho} P^*_T f) = \Lambda_{\rho,T}f$ and $\Lambda^\ast_{\rho,T} = R \Lambda_{\rho,T} R$ on $L^2((0,T)\times\partial\Omega)$, we have
\begin{align*}
I(T,T) & = (f,J \Lambda_{\rho} P^\ast_T h)_{L^2((0,T)\times\partial\Omega)} - (\Lambda_{\rho,T} f, J P^\ast_T h)_{L^2((0,T)\times\partial\Omega)} \\
 & = (f,J \Lambda_{\rho} P^\ast_T h)_{L^2((0,T)\times\partial\Omega)} - (f, R \Lambda_{\rho,T} R J P^\ast_T h)_{L^2((0,T)\times\partial\Omega)} \\
 & = (f,Kh)_{L^2((0,T)\times\partial\Omega)}
\end{align*}

For general $f, h \in L^2((0,T)\times\partial\Omega)$, notice that $K: L^2((0,T)\times\partial\Omega) \rightarrow L^2((0,T)\times\partial\Omega)$ is a continuous operator since all the operators in the definition~\eqref{eq:K} are continuous. The result follows from the density of compactly supported smooth functions in $L^2$.
\end{proof}

\begin{remark} \label{rm:extension}
The zero extension $P^*_T$ plays no special role than other extensions in Proposition~\ref{thm:waveinner}. In fact, if $Q^*_T$ is another extension such that $P_T Q^*_T f = f$ for all $f \in L^2((0,T)\times\partial\Omega)$, then 
$$
u^{P^*_T f} = u^{Q^*_T f} \quad \text{ on } [0,T]\times\Omega
$$ 
as they satisfy identical initial and boundary conditions on $[0,T]$. Because of this observation, we often omit the notation for extension and simply write $u^{P^*_T f}$ as $u^f$ and write $\dot{\Lambda}_T$ as $\dot\Lambda$ when considering functions on $[0,T]$.
\end{remark}

\bigskip

\begin{coro}
Suppose $f, h\in C^\infty_c((0,T]\times\partial\Omega)$. Then
\begin{equation} \label{eq:ufuhpp}
( (\Delta-q) u^{f}(T), u^h(T))_{L^2(\Omega)} = (\partial^2_t f,Kh)_{L^2((0,T)\times\partial\Omega)} = (K \partial^2_t f,h)_{L^2((0,T)\times\partial\Omega)}.
\end{equation}
\end{coro}
\begin{proof}
Replacing $f$ by $\partial^2_t f$ in~\eqref{eq:ufuh}, we get
$$
(\rho u^{\partial^2_t f}(T), u^h(T))_{L^2(\Omega)} = (\partial^2_t f,Kh)_{L^2((0,T)\times\partial\Omega)} = (K\partial^2_t f,h)_{L^2((0,T)\times\partial\Omega)}.
$$
As both $u^{\partial^2_t f}$ and $\partial^2_t u^f$ satisfy~\eqref{eq:bvp} with $f$ replaced by $\partial^2_t f$, they must be equal thanks to the well-posedness of the boundary value problem. We conclude
$$
\rho u^{\partial^2_t f} = \rho \partial^2_t u^f
= \Delta u^f - q u^f.
$$
\end{proof}


Recall that $\Lambda_{\rho}=\Lambda_{\rho_0} + \epsilon \dot{\Lambda}_{\dot{\rho}}$ in the linearization setting. When there is no risk of confusion, we write  $\dot{\Lambda}_{\dot\rho}$ simply as $\dot{\Lambda}$.

Accordingly, we decompose $K=K_0+\epsilon\dot{K}$. Here
$K_0$ is the connecting operator for the background medium: 
\begin{equation} \label{eq:K0}
K_0:= J \Lambda_{\rho_0} P^\ast_T - R \Lambda_{\rho_0,T} R J P^\ast_T.
\end{equation}
$K_0$ can be explicitly computed since $\Lambda_{\rho_0}$ is known.
$\dot{K}$ is the resulting perturbation in the connecting operator:
\begin{equation} \label{eq:Kdot}
\dot{K}:= J \dot{\Lambda} P^\ast_T - R \dot{\Lambda}_{T} R J P^\ast_T.
\end{equation}
where $\dot\Lambda_{T}f := P_T \dot\Lambda P^*_T f$. Note that $\dot{K}$ can be explicitly computed once $\dot{\Lambda}$ is given.

Let us introduce some function spaces in order to discuss the mapping properties of $\dot{\Lambda}$ and $\dot K$. Denote 
\begin{align*}
    H^1_{cc}((0,T)\times\partial\Omega) := & \{f\in H^1((0,T)\times\partial\Omega): f(0,x) = 0 \text{ for all } x\in\partial\Omega\} \\
    H^2_{cc}((0,T)\times\partial\Omega) := & \{f\in H^2((0,T)\times\partial\Omega): f(0,x) = \partial_t f(0,x) = 0 \text{ for all } x\in\partial\Omega\}.
\end{align*}
and equip them with the usual $H^1$-norm and $H^2$-norm, respectively.

\begin{lemma} \label{thm:Lambdadotspace}
    \begin{align*}
    \dot\Lambda_T: & \; H^1_{cc}((0,T)\times\partial\Omega) \to L^2((0,T)\times\partial\Omega)  \\
    \dot\Lambda_T: & \; H^2_{cc}((0,T)\times\partial\Omega) \to H^1((0,T)\times\partial\Omega)
    \end{align*}
    are bounded linear operators.    
\end{lemma}

\begin{proof}
    By~\cite[Theorem A(2)(4)]{lasiecka1991regularity}, we have that
    \begin{align*}
    \dot\Lambda: & \; H^1_{cc}((0,2T)\times\partial\Omega) \to L^2((0,2T)\times\partial\Omega)  \\
    \dot\Lambda: & \; H^2_{cc}((0,2T)\times\partial\Omega) \to H^1((0,2T)\times\partial\Omega)
    \end{align*}
    are bounded linear operators.
    By Remark~\ref{rm:extension}, we can choose the extension $P^*_T$ to be regularity-preserving in the sense that
    \begin{align*}
            P^*_T: & \; H^1_{cc}((0,T)\times\partial\Omega) \rightarrow H^1_{cc}((0,2T)\times\partial\Omega) \\           P^*_T: & \; H^2_{cc}((0,T)\times\partial\Omega) \rightarrow H^2_{cc}((0,2T)\times\partial\Omega).
    \end{align*}
    Therefore, $\dot{\Lambda}_T = P_T \dot\Lambda P^*_T$ has the desired mapping properties.
\end{proof}

\begin{lemma} \label{thm:Jspace}
    \begin{align*}
        J: & \; L^2([0,2T]\times\partial\Omega) \rightarrow L^2([0,T]\times\partial\Omega) \\
        J: & \; H^1_{cc}((0,2T)\times\partial\Omega) \rightarrow H^1((0,T)\times\partial\Omega)
    \end{align*}
    are bounded linear operators.
\end{lemma}
\begin{proof}
    The first claim is a consequence of the following estimate
    \begin{align*}                      \|Jf\|_{L^2([0,T]\times\partial\Omega)}^2=&\int_0^T\int_{\partial\Omega}\left|\frac{1}{2}\int_t^{2T-t}f(\tau,x)\dif\tau\right|^2\dif x\dif t\\
    =&\frac{1}{4}\int_0^T\int_{\partial\Omega}\left(\int_t^{2T-t}\left|f(\tau,x)\right|\dif\tau\right)^2\dif x\dif t\\
    \leq&\frac{1}{4}\int_0^T\int_{\partial\Omega} 2(T-t) \int_t^{2T-t}\left|f(\tau,x)\right|^2\dif\tau\dif x\dif t\\
    \leq & \frac{1}{4}\int_0^T\int_{\partial\Omega} 2T \int_0^{2T}\left|f(\tau,x)\right|^2\dif\tau\dif x\dif t \\
    =&\frac{T^2}{2}\|f\|_{L^2([0,2T]\times\partial\Omega)}^2.
    \end{align*}

    To prove the second claim, notice that $\partial_{x_j} Jf = J (\partial_{x_j}f)$ for all $j=1,2,\dots,n$. Hence the $L^2$-estimate above gives
    $$
    \|\partial_{x_j} Jf\|^2_{L^2((0,T)\times\partial\Omega)} = \|J ( \partial_{x_j} f)\|^2_{L^2((0,T)\times\partial\Omega)} \leq \frac{T^2}{2} \|\partial_{x_j} f\|_{L^2((0,2T)\times\partial\Omega)}^2.
    $$
    On the other hand, using $f(0,x)=0$ for $f\in H^1_{cc}((0,2T)\times\partial\Omega)$, we have
    \begin{align*}
        \partial_t Jf(t) & = -\frac{1}{2} \left[ f(2T-t) + f(t) \right] = -\frac{1}{2} \left[ \int^{2T-t}_0 \partial_\tau f(\tau,x) \dif \tau + \int^{t}_0 \partial_\tau f(\tau,x) \dif \tau \right] \\
        = & -J(\partial_t f)(t) - \int^t_0 \partial_\tau f(\tau,x) \dif \tau, \qquad t\in (0,T).
    \end{align*}
Therefore, 
    \begin{align*}
        \|\partial_t Jf\|^2_{L^2((0,T)\times\partial\Omega)} \leq & 2 \|J (\partial_t f)\|^2_{L^2((0,T)\times\partial\Omega)} + 
        2 \int_0^T\int_{\partial\Omega}\left|\int^t_{0} \partial_\tau f(\tau,x)\dif\tau\right|^2\dif x\dif t\\
        \leq & T^2 \|\partial_t f\|^2_{L^2((0,2T)\times\partial\Omega)} + 
        2 \int_0^T\int_{\partial\Omega} t \int^t_{0} |\partial_\tau f(\tau,x)|^2 \dif\tau \dif x\dif t\\
        \leq & T^2 \|\partial_t f\|^2_{L^2((0,2T)\times\partial\Omega)} + 
        2 \int_0^T\int_{\partial\Omega} t \int^{2T}_{0} |\partial_\tau f(\tau,x)|^2 \dif\tau \dif x\dif t\\
        = & T^2 \|\partial_t f\|^2_{L^2((0,2T)\times\partial\Omega)} + T^2 \|\partial_t f\|^2_{L^2((0,2T)\times\partial\Omega)} \\
        = & 2 T^2 \|\partial_t f\|^2_{L^2((0,2T)\times\partial\Omega)}.
    \end{align*}
This completes the proof.
\end{proof}

\begin{lemma} \label{thm:Kdotspace}
    $$
    \dot K: H^1_{cc}((0,T)\times\partial\Omega) \to L^2((0,T)\times\partial\Omega)
    $$
    is a bounded linear operator.
\end{lemma}

\begin{proof}
Based on the definition of $\dot K$ in~\eqref{eq:Kdot}, we analyze the two terms $J\dot\Lambda P^*_T$ and $R\dot\Lambda_T R J P^*_T$ separately.

By Remark~\ref{rm:extension}, we can choose the extension $P^*_T$ to be regularity-preserving so that $P^*_T: H^1_{cc}((0,T)\times\partial\Omega) \rightarrow H^1_{cc}((0,2T)\times\partial\Omega)$ is bounded.
By~\cite[Theorem A(2)(4)]{lasiecka1991regularity}, $\dot\Lambda: H^1_{cc}((0,2T)\times\partial\Omega) \to L^2((0,2T)\times\partial\Omega)$ is bounded.
By Lemma~\ref{thm:Lambdadotspace}, $J: L^2((0,2T)\times\partial\Omega) \to L^2((0,T)\times\partial\Omega)$ is bounded.
Hence $J\dot\Lambda P^*_T: H^1_{cc}((0,T)\times\partial\Omega) \rightarrow L^2((0,T)\times\partial\Omega)$ is bounded.

For the second term, let us take an arbitrary $f\in H^1_{cc}((0,T)\times\partial\Omega)$. $P^*_T$ preserves the $H^1$-regularity as we have seen. $R$ is clearly an isometry on $H^1((0,T)\times\partial\Omega)$. $J$ preserves the $H^1$-regularity by Lemma~\ref{thm:Jspace}. Moreover, $RJP^*_Tf(0,x) = JP^*_Tf(T,x) = 0$ for all $x\in\partial\Omega$. Hence $RJP^*_T f \in H^1_{cc}((0,T)\times\partial\Omega)$. By Lemma~\ref{thm:Lambdadotspace}, $\dot\Lambda_T: H^1_{cc}((0,T)\times\partial\Omega) \to L^2((0,T)\times\partial\Omega)$ is bounded. $R$ is also an isometry on $L^2((0,T)\times\partial\Omega)$. Therefore, the composition $R\dot\Lambda_T R J P^*_T: H^1_{cc}((0,T)\times\partial\Omega) \rightarrow L^2((0,T)\times\partial\Omega)$ is bounded.
\end{proof}

\section{Integral Identity and Controllability} \label{sec:id}


We derive an integral identity in Proposition~\ref{thm:id} that is essential for our linearized BC method. This identity can be understood as the linearized Blagove\u{s}\u{c}enski\u{ı}’s identity with a free parameter.

\begin{prop} \label{thm:id}
Let $0\neq \lambda\in\mathbb{R}$ be a nonzero real number. If $f, h\in C^\infty_c((0,T]\times\partial\Omega)$ satisfy
\begin{equation} \label{eq:fhvanish}
[\Delta - q + \lambda \rho_0 ] u_0^{f}(T) =
[\Delta - q + \lambda \rho_0 ] u_0^h(T) = 0
\quad \text{ in } \Omega,
\end{equation}
then the following identity holds:
\begin{equation} \label{eq:keyid}
-\lambda (\dot{\rho} u_0^f(T), u_0^h(T))_{L^2(\Omega)} = (\partial^2_t f + \lambda f,\dot{K}h)_{L^2((0,T)\times\partial\Omega)} + (\dot{\Lambda} f(T), h(T))_{L^2(\partial\Omega)}.
\end{equation}
\end{prop}

\begin{proof}

For $f, h\in C^\infty_c((0,T]\times\partial\Omega)$, we will make use of~\eqref{eq:ufuh} ~\eqref{eq:ufuhpp} to obtain some identities. First, let $\epsilon=0$ in~\eqref{eq:ufuh} we obtain
$$
(\rho_0 u_0^f(T), u_0^h(T))_{L^2(\Omega)} = (f,K_0h)_{L^2((0,T)\times\partial\Omega)} = (K_0f,h)_{L^2((0,T)\times\partial\Omega)}.
$$
Next, differentiating~\eqref{eq:ufuh} in $\epsilon$ and let $\epsilon=0$, we obtain
\begin{align}
 & (f,\dot{K}h)_{L^2((0,T)\times\partial\Omega)} = (\dot{K}f,h)_{L^2((0,T)\times\partial\Omega)} \nonumber \\
= & (\dot{\rho} u_0^f(T), u_0^h(T))_{L^2(\Omega)} + (\rho_0 \dot{u}^f(T), u_0^h(T))_{L^2(\Omega)} + (\rho_0 u_0^f(T), \dot{u}^h(T))_{L^2(\Omega)}  \label{eq:ufuhlinearized}.
\end{align}

Similarly, let $\epsilon=0$ in~\eqref{eq:ufuhpp} we obtain
$$
(\Delta u_0^{f}(T) - q u_0^{f}(T), u_0^h(T))_{L^2(\Omega)} = (\partial^2_t f,K_0h)_{L^2((0,T)\times\partial\Omega)} = (K_0 \partial^2_t f,h)_{L^2((0,T)\times\partial\Omega)}.
$$
Next, differentiating~\eqref{eq:ufuhpp} in $\epsilon$ and letting $\epsilon=0$, we obtain
\begin{align*}
 & (\partial^2_t f,\dot{K}h)_{L^2((0,T)\times\partial\Omega)} = (\dot{K} \partial^2_t f,h)_{L^2((0,T)\times\partial\Omega)} \\
 = & ( (\Delta-q) \dot{u}^{f}(T), u_0^h(T))_{L^2(\Omega)}  +  ( (\Delta-q) u_0^{f}(T), \dot{u}^h(T))_{L^2(\Omega)} \\
 = & ( \dot{u}^{f}(T), (\Delta-q)  u_0^h(T))_{L^2(\Omega)} - ( \dot{\Lambda}f(T), h(T))_{L^2(\partial\Omega)} + ( (\Delta-q) u_0^{f}(T), \dot{u}^h(T))_{L^2(\Omega)}
\end{align*}
where the last inequality follows from integration by parts along with the fact that $\partial_\nu \dot{u}=0$ and $\dot{u}^{f}|_{(0,2T)\times\partial\Omega} = \dot{\Lambda} f$, and $L^2(\Omega) = L^2(\Omega,dx)$ is the $L^2$-space equipped with the usual Lebesgue measure. Add~\eqref{eq:ufuhlinearized} multiplied by $\lambda\in\mathbb{R}$ to get
\begin{align*}
 & (\partial^2_t f + \lambda f, \dot{K}h)_{L^2((0,T)\times\partial\Omega)} + (\dot{\Lambda} f(T), h(T))_{L^2(\partial\Omega)} \nonumber \\
= & (\partial^2_t f , \dot{K}h)_{L^2((0,T)\times\partial\Omega)} + (\dot{\Lambda} f(T), h(T))_{L^2(\partial\Omega)} + ( \lambda f, \dot{K}h)_{L^2((0,T)\times\partial\Omega)} \\
= & (\dot{u}^{f}(T), [\Delta - q + \lambda \rho_0 ] u_0^h(T))_{L^2(\Omega)} + ([\Delta - q + \lambda \rho_0 ] u^{f}_0(T),  \dot{u}^h(T))_{L^2(\Omega)} + \lambda (\dot{\rho} u_0^f(T), u_0^h(T))_{L^2(\Omega)}.
\end{align*}
If $[\Delta - q + \lambda \rho_0 ] u_0^{f}(T) =
[\Delta - q + \lambda \rho_0 ] u_0^h(T) = 0
\text{ in } \Omega$, the first term and second term on the right-hand side vanish, resulting in~\eqref{eq:keyid}.
\end{proof}

\medskip
The following boundary control estimate is established in~\cite{oksanen2022linearized}. 
Given a strictly positive $\rho_0\in C^\infty(\overline{\Omega})$, we will write $g := \rho_0 dx^2$ for the Riemannian metric associated to $\rho_0$, and denote by $S\overline{\Omega}$ the unit sphere bundle over the closure $\overline{\Omega}$ of $\Omega$.

\begin{prop}[\cite{oksanen2022linearized}] \label{thm:control}
Let $\rho_0\in C^\infty(\overline{\Omega})$ be strictly positive and $q_0\in C^\infty(\overline{\Omega})$.
Suppose that all maximal\footnote{For a maximal geodesic $\gamma : [a, b] \to \overline\Omega$ there may exists $t \in (a,b)$ such that $\gamma(t) \in \p \Omega$. The geodesics are maximal on the closed set $\overline{\Omega}$.} geodesics on $(\overline{\Omega}, g)$ have length strictly less than some fixed $T > 0$.  
Then for any $\phi \in C^\infty(\overline{\Omega})$ there is $f \in C_c^\infty((0,T] \times \p \Omega)$ such that 
    \begin{equation}\label{eq:controleqn}
u_0^f(T) = \phi \quad \text{in $\Omega$},
    \end{equation}
where $u_0$ is the solution of~\eqref{eq:unperturbedBVP}.
Moreover, there is $C>0$, independent of $\phi$, such that
    \begin{equation}\label{control_estimate}
\norm{f}_{H^2((0,T) \times \p \Omega)} 
\le C \norm{\phi}_{H^4(\Omega)}.
    \end{equation}
\end{prop}

\section{Stability and Reconstruction} \label{sec:stabandrecon}

We derive the stability estimate and reconstruction procedure in this section.

\subsection{Case 1: $\rho_0=$ const} \label{sec:q0vanish}


\bigskip
Without loss of generality, we assume $\rho_0=1$. This case can be coped with using a similar method as in~\cite{oksanen2022linearized}. Here, we simply sketch the idea for dimension $n\geq 2$. Details can be found in~\cite[Section 4]{oksanen2022linearized}.

Let $\lambda>0$. The equations~\eqref{eq:fhvanish} become the perturbed Helmholtz equation
$$
[\Delta - q + \lambda ] u_0^{f}(T) =
[\Delta - q + \lambda ] u_0^h(T) = 0 \quad \text{ in } \Omega.
$$
We separate the discussion for $q \equiv 0$ and $q\not\equiv 0$.

When $q\equiv 0$, let $\theta\in\mathbb{S}^{n-1}$ be an arbitrary unit vector. Thanks to Proposition~\ref{thm:control}, the boundary control equations
$$
u^f_0(T) = u^h_0(T) = e^{i\sqrt{\lambda}\theta\cdot x}
$$
admit solutions $f,h\in C^\infty_c((0,T]\times\partial\Omega)$.
Inserting such $f,h$ into~\eqref{eq:keyid} yields $\mathcal{F}[\dot \rho](2\sqrt{\lambda}\theta)$, the Fourier transform of $\dot\rho$ at $2\sqrt{\lambda}\theta$. Varying $\lambda>0$ and $\theta\in\mathbb{S}^{n-1}$ recovers the full Fourier transform of $\dot{\rho}$. The procedure is summarized in Algorithm~\ref{alg:constq}.

\begin{algorithm}[h] 
    \SetAlgoLined   
	\caption{Linearized Boundary Control Reconstruction of $\dot{\rho}$ when $q\equiv 0$}  \label{alg:constq}
 \vspace{1ex}
	\KwIn{low-pass filter $J$, time-reversal operator $R$, projection operator $P_T$, linearized ND map $\dot{\Lambda}_{\dot{\rho}}$ \medskip}
	
		 Choose $\lambda>0$ and $\theta\in\mathbb{S}^{n-1}$.
		 
	   Solve the boundary control equations $u^f_0(T) = u^h_0(T) = e^{i\sqrt{\lambda}\theta\cdot x}$ for $f$ and $h$. 

          Construct the linearized connecting operator $\dot{K}$ by $\dot{K}:= J \dot{\Lambda} P^\ast_T - R \dot{\Lambda}_{T} R J P^\ast_T$.

		 Compute $\mathcal{F}[\dot \rho](2\sqrt{\lambda}\theta)$ by
\begin{equation} \label{eq:rhodotrecon1}
\mathcal{F}[\dot \rho](2\sqrt{\lambda}\theta) = -\frac{1}{\lambda} \left[
(\partial^2_t f + \lambda f,\dot{K}h)_{L^2((0,T)\times\partial\Omega)} + (\dot{\Lambda} f(T), h(T))_{L^2(\partial\Omega)}
\right].
\end{equation}

         Repeat the above steps with various $\lambda>0$ and $\theta\in\mathbb{S}^{n-1}$ to recover the Fourier transform $\mathcal{F}[\dot \rho]$.
         
	     Invert the Fourier transform to recover $\dot\rho$. \medskip

    	\KwOut{sound speed perturbation $\dot{\rho}$}

\end{algorithm}

\begin{remark}
Step 2 of Algorithm~\ref{alg:constq} requires solving the boundary control equations, for which a solution exists due to Proposition~\ref{thm:control}. As $\rho_0\equiv 1$ and $q$ is known, we can solve these linear systems to find $f$ and $h$. In Section~\ref{sec:TR}, we describe a time-reversal procedure to solve these equations in a special 1D setup.
\end{remark}

\bigskip
The following Lipschitz-type stability estimate can be readily derived from the reconstruction formula~\eqref{eq:rhodotrecon1}.  The proof is nearly a word-by-word repetition of~\cite[Theorem 7]{oksanen2022linearized} and is omitted.

\begin{thm} \label{thm:stab1}
Suppose $q\equiv 0$. There exists a constant $C>0$, independent of $\lambda$, such that 
$$
\left| \mathcal{F}[{\dot{\rho}}](\sqrt{2\lambda}\theta) \right| \leq C (2+\lambda)\lambda^3  
\left(
\|\dot{K}\|_{H^1_{cc}((0,T)\times\partial\Omega)\rightarrow L^2((0,T)\times\partial\Omega)} + \|\dot{\Lambda}\|_{H^2_{cc}((0,T)\times\partial\Omega) \rightarrow H^1((0,T)\times\partial\Omega)}
   \right)
$$
Here $\dot{K}$ is viewed as a linear function of $\dot{\Lambda}$ as is defined in~\eqref{eq:Kdot}. \label{thm:nopotential}
\end{thm}

\bigskip
When $q\not\equiv 0$, 
let $\theta,\omega\in\mathbb{R}^n$ be two vectors such that $\theta\perp\omega$, and let $k,l>0$ be two real numbers such that $k^2+l^2=\lambda$. Define two functions 
\begin{align*}
    \phi_+(x) := e^{i(k\theta+l\omega)\cdot x} + r_+(x;\lambda), \quad\quad
    \phi_-(x) := e^{i(k\theta-l\omega)\cdot x} + r_-(x;\lambda)
\end{align*}
where $r_\pm$ solve the equation $(\Delta+\lambda-q) r_\pm = q e^{i(k\theta\pm l\omega)\cdot x}$ and the outgoing Sommerfeld radiation condition. 
It is shown~\cite[Lemma 13]{oksanen2022linearized} that such solutions $r_\pm$ exist, and their Sobolev norms of order $s>0$ have the following asymptotic behavior
$$
\|r_\pm\|_{H^s(\mathbb{R}^n)} = O(\lambda^{\frac{s-1}{2}}), \quad\quad \text{ as } \lambda\rightarrow \infty.
$$
Thanks to Proposition~\ref{thm:control}, the boundary control equations
$$
u^f_0(T) = \phi_+, \quad\quad 
u^h_0(T) = \phi_-
$$
admit solutions $f,h\in C^\infty_c((0,T]\times\partial\Omega)$.
Inserting such $f,h$ into~\eqref{eq:keyid} then taking the limit $l \rightarrow\infty$ yields $\mathcal{F}[\dot\rho](2k\theta)$. Varying $k>0$ and $\theta\in\mathbb{S}^{n-1}$ recovers the full Fourier transform of $\dot\rho$. The procedure is summarized in Algorithm~\ref{alg:varq}.

\begin{algorithm}[h]
    \SetAlgoLined
	\caption{Linearized Boundary Control Reconstruction of $\dot{\rho}$ when $n\geq 2$ and $q\not\equiv 0$}  \label{alg:varq}
 \vspace{1ex}
	\KwIn{low-pass filter $J$, time-reversal operator $R$, projection operator $P_T$, linearized ND map $\dot{\Lambda}_{\dot{\rho}}$ \medskip}
	
		 Choose $k,l>0$ with $k^2+l^2=\lambda^2$ and $\theta,\omega \in\mathbb{S}^{n-1}$ with $\theta\perp\omega$.

          Solve the equation $(\Delta+\lambda-q) r_\pm = q e^{i(k\theta\pm l\omega)\cdot x}$ along with the outgoing Sommerfeld radiation condition to find $r_\pm$.
   
	   Solve the boundary control equation 
    $$
    u^f_0(T) = e^{i(k\theta+l\omega)\cdot x} + r_+(x;\lambda), \quad\quad
    u^h_0(T) = e^{i(k\theta-l\omega)\cdot x} + r_-(x;\lambda)
    $$
        for $f$ and $h$. 

        Construct the linearized connecting operator $\dot{K}$ by
$\dot{K}:= J \dot{\Lambda} P^\ast_T - R \dot{\Lambda}_{T} R J P^\ast_T$.

		 Compute $\mathcal{F}[\dot \rho](2k \theta)$ by
\begin{equation} \label{eq:rhodotrecon2}
\mathcal{F}[\dot \rho](2 k \theta) = - \lim_{l\rightarrow\infty} \frac{1}{k^2+l^2} \left[
(\partial^2_t f + (k^2+l^2)) f,\dot{K}h)_{L^2((0,T)\times\partial\Omega)} + (\dot{\Lambda} f(T), h(T))_{L^2(\partial\Omega)}
\right].
\end{equation}

         Repeat the above steps with various $\lambda>0$ and $\theta\in\mathbb{S}^{n-1}$ to recover the Fourier transform $\mathcal{F}[\dot \rho]$.
         
	     Invert the Fourier transform to recover $\dot\rho$. \medskip

    	\KwOut{sound speed perturbation $\dot{\rho}$}

\end{algorithm}

\subsection{Case 2: $\rho_0 \neq \text{constant}$.}


When $\rho_0=\rho_0(x)>0$ is non-constant, the equations~\eqref{eq:fhvanish} are no longer perturbed Helmholtz equations, but Schr\"odinger's equations with the potential $-q+\lambda \rho_0 \in L^\infty(\Omega)$. 
The idea is to employ Schr\"odinger solutions to probe  based on the identity~\eqref{eq:keyid}.

The class of solutions we will resort to are the \textit{complex geometric optics (CGO) solutions} that were first proposed in~\cite{sylvester1987global} for dimension $n\geq 3$. A CGO solution $\phi$ is a function of the form
\begin{equation} \label{eq:CGO}
    \phi(x) := e^{i\zeta\cdot x}(1+r(x)).
\end{equation}
where $\zeta\in\mathbb{C}^n$ is a complex vector with $\zeta\cdot\zeta=0$, and the remainder term $r(x)$ satisfies
$$
    \Delta r + 2\zeta\cdot\nabla r - (q - \lambda\rho_0) r = q - \lambda\rho_0.
$$
Moreover, $r\rightarrow 0$ in a certain function space as $|\zeta|\rightarrow\infty$.

The following proposition is a direct application of~\cite[Theorem 2.3 and Corollary 2.4]{sylvester1987global} to the Schr\"odinger's equation $(\Delta - q + \lambda\rho_0) \phi=0$.

\begin{lemma}[{\cite[Theorem 2.3 and Corollary 2.4]{sylvester1987global}}]\label{thm:CGOregularity}
Let $n\geq 3$ and $s\in\mathbb{R}$ a real number such that $s>\frac{n}{2}$.
Let $\zeta\in\mathbb{C}^n$ be a complex vector with $\zeta\cdot\zeta=0$ and $|\zeta|\geq\varepsilon_0>0$ for some positive constant $\varepsilon_0$ 
 There exist positive constants $C_0$, $C_1$, depending on $s,n,\varepsilon_0$ and $\Omega$, such that if 
\[C_0\|q - \lambda\rho_0\|_{H^{s}(\Omega)}<|\zeta|,\]
then $\phi=\phi(x)$ defined in~\eqref{eq:CGO}
satisfies $(\Delta - q + \lambda\rho_0) \phi=0$; moreover
$$
    \|r\|_{H^{s}(\Omega)}\leq \frac{C_1}{|\zeta|}\|q - \lambda\rho_0\|_{H^{s}(\Omega)}.
$$
\end{lemma}

We now construct specific CGO solutions that are useful for our purpose. Let $\xi\in\mathbb{R}^{n} (n\geq 3)$ be an arbitrary non-zero vector, and let $e^{(1)},e^{(2)}\in\mathbb{S}^{n-1}$ be two real unit vectors such that $\{\xi,e^{(1)},e^{(2)} \}$ forms an orthogonal set. Choose a positive number $R$ with $R\geq \frac{|\xi|}{\sqrt{2}}$. Define
$$
    \zeta^{(1)} :=-\frac{1}{2}\xi + i\frac{R}{\sqrt{2}}e^{(1)} + \sqrt{\frac{R^2}{2}-\frac{|\xi|^2}{4}}e^{(2)}, \quad\quad 
    \zeta^{(2)} :=-\frac{1}{2}\xi - i\frac{R}{\sqrt{2}}e^{(1)} - \sqrt{\frac{R^2}{2}-\frac{|\xi|^2}{4}}e^{(2)}.
$$
It is easy to verify that
\[\zeta^{(1)}+\zeta^{(2)}=-\xi,\quad\zeta^{(j)}\cdot\zeta^{(j)}=0,\quad|\zeta^{(j)}|=R,\quad\textup{for }j=1,2.\]
If $R$ is sufficiently large, by Lemma~\ref{thm:CGOregularity}, we can construct CGO solutions
\begin{equation} 
\phi_j(x) = e^{i\zeta^{(j)}\cdot x}(1+r_j(x))
\end{equation}
where the remainder term $r_j$ satisfies
\begin{equation} \label{eq:rbound}
\|r_j\|_{H^{s}(\Omega)}\leq \frac{C_1}{|\zeta^{(j)}|}\|q - \lambda\rho_0\|_{H^{s}(\Omega)}\leq \frac{C_1}{C_0}. 
\end{equation}
(Here, $C_0$ is the constant in Lemma~\ref{thm:CGOregularity}.)

Thus for $s>\frac{n}{2}$,
\[\|\phi_j\|_{H^{s}(\Omega)}\leq\|e^{i\zeta^{(j)}\cdot x}\|_{H^{s}(\Omega)}\|1+r_j\|_{H^{s}(\Omega)}\leq\left(|\Omega|^{\frac{1}{2}}+\frac{C_1}{C_0}\right)\|e^{i\zeta^{(j)}\cdot x}\|_{H^{s}(\Omega)}.\]
By choosing $\lambda_0$ such that for any $\lambda>\lambda_0$, we have \[C_0\|q - \lambda\rho_0\|_{H^{s}(\Omega)}\geq\frac{1}{\sqrt{n}},\]
\begin{align*}
    \|e^{i\zeta^{(j)}\cdot x}\|_{H^k(\Omega)}^2=&\sum_{i=0}^k\sum_{|\alpha|=i}\|D^\alpha e^{i\zeta^{(j)}\cdot x}\|_{L^2(\Omega)}^2=\sum_{i=0}^k\sum_{\sum_{m=1}^n\alpha_m=i}\left(\prod_{m=1}^n|\zeta^{(j)}_m|^{\alpha_m}\right)\| e^{\textup{Im}\zeta^{(j)}\cdot x}\|_{L^2(\Omega)}^2\\
    \leq&C\sum_{i=0}^k\|\zeta^{(j)}\|_1^ie^{\sqrt{2}R}\leq C\sum_{i=0}^k(\sqrt{n}\|\zeta^{(j)}\|_2)^ie^{\sqrt{2}R}=C\sum_{i=0}^k(\sqrt{n}R)^ie^{\sqrt{2}R}\leq CR^ke^{\sqrt{2}R}
\end{align*}
where $C$ only depend on $k,n,\Omega$. Using an interpolation argument, we obtain
\begin{equation}\label{eq:CGOHsnorm}
    \|\phi_j\|_{H^{s}(\Omega)}\leq\left(|\Omega|^{\frac{1}{2}}+\frac{C_1}{C_0}\right)\|e^{i\zeta^{(j)}\cdot x}\|_{H^{s}(\Omega)}\leq CR^{\frac{s}{2}}e^{\frac{R}{\sqrt{2}}}.
\end{equation}

\bigskip
We are ready to derive some stability estimates. 
For simplicity we denote 
\[\delta :=   \|\dot{K}\|_{H^1_{cc}((0,T)\times\partial\Omega)\rightarrow L^2((0,T)\times\partial\Omega)} + \|\dot{\Lambda}\|_{H^2_{cc}((0,T)\times\partial\Omega) \rightarrow H^1((0,T)\times\partial\Omega)}.
   \]
These norms are valid in view of Lemma~\ref{thm:Lambdadotspace} and Lemma~\ref{thm:Kdotspace}.
We begin with a pointwise estimate for $\dot{\rho}$ in the Fourier domain.

\medskip
\begin{lemma} \label{thm:fourierest}
    Let $s>\frac{n}{2}$ with $n\geq 3$. Suppose there exists a constant $M>0$ such that
    $$
    \|q\|_{H^{\max(s,4)}(\Omega)} \leq M, \quad \|\rho_0\|_{H^{\max(s,4)}(\Omega)} \leq M.
    $$
    Then there exists a constant $C$, independent of $\lambda$ and $\delta$, such that
    \begin{equation}\label{eq:fourierest}
    |\hat{\dot\rho}(\xi)|\leq
        \begin{cases} 
            C \left[ \frac{\lambda+1}{\lambda}(\lambda+1)^{\max(s,4)}e^{\sqrt{2}a_0(\lambda+1)}\delta+\frac{1}{a_0}\|\dot{\rho}\|_{H^{-s}(\Omega)} \right] &|\xi|\leq\sqrt{2}a_0(\lambda+1) 
 \vspace{1ex} \\
            C \left[\frac{\lambda+1}{\lambda}|\xi|^{\max(s,4)}e^{|\xi|}\delta+\frac{\lambda+1}{|\xi|}\|\dot{\rho}\|_{H^{-s}(\Omega)} \right] &|\xi|\geq\sqrt{2}a_0(\lambda+1)
        \end{cases}
    \end{equation}
    for any $\lambda>0$ and sufficiently small $\delta$. Here, $a_0$ is a constant satisfying $a_0 \geq C_0M$, where $C_0$ is the constant in Lemma~\ref{thm:CGOregularity}. 
\end{lemma}
\begin{proof}

From Proposition~\ref{thm:control}, there exist boundary controls $f_j$ such that $u_0^{f_j}(T)=\phi_j$ for the CGO solutions $\phi_j$ defined in~\eqref{eq:CGO}. 
By Proposition~\ref{thm:id}, we have
\[
    \begin{aligned}
    &\left|\int_\Omega\dot{\rho}\phi_1\phi_2\dif x\right| \vspace{1ex}\\
    = & \frac{1}{\lambda}\left|(\partial^2_t f_1 + \lambda f_1,\dot{K}f_2)_{L^2((0,T)\times\partial\Omega)} + (\dot{\Lambda} f_1(T), f_2(T))_{L^2(\partial\Omega)}\right| \vspace{1ex}\\
     \leq & \frac{1}{\lambda}\left[\|\partial^2_t f_1 + \lambda f_1\|_{L^2((0,T)\times\partial\Omega)}\|\dot{K}f_2\|_{L^2((0,T)\times\partial\Omega)} + 
    \|\dot{\Lambda}f_1(T)\|_{L^2(\partial\Omega)} \|f_2(T)\|_{L^2(\partial\Omega)} \right] \\  
    \leq & \frac{1}{\lambda}\big[ (1+\lambda)\|f_1\|_{H^2((0,T)\times\partial\Omega)}\|\dot{K}f_2\|_{L^2((0,T)\times\partial\Omega)} + \|\dot{\Lambda} f_1\|_{H^1((0,T)\times\partial\Omega)}\| f_2\|_{H^1((0,T)\times\partial\Omega)}\big]\\
    \leq & \frac{1}{\lambda}\big[ (1+\lambda)\|f_1\|_{H^2((0,T)\times\partial\Omega)} \|\dot{K}\|_{H^1_{cc}((0,T)\times\partial\Omega)\rightarrow L^2((0,T)\times\partial\Omega)} \|f_2\|_{H^1((0,T)\times\partial\Omega)}   \\
    & + 
    \|\dot{\Lambda} \|_{H^2_{cc}((0,T)\times\partial\Omega)\rightarrow H^1((0,T)\times\partial\Omega)}
\|f_1\|_{H^2_{}((0,T)\times\partial\Omega)}\| f_2\|_{H^2((0,T)\times\partial\Omega)}\big]\\
    \leq&C_\lambda \delta   
    \end{aligned}
\]
where the first inequality is by the Cauchy-Schwarz inequality, the second inequality by the trace estimate, and the last inequality by Lemma~\ref{thm:Kdotspace}.
Here, the constant $C_\lambda$ is
\begin{equation}
    \begin{aligned}
        C_\lambda :=& \frac{1+\lambda}{\lambda}\|f_1 \|_{H^2((0,T)\times\partial\Omega)}
    \|f_2 \|_{H^2((0,T)\times\partial\Omega)} \\
    \leq& C \frac{1+\lambda}{\lambda} \|\phi_1\|_{H^4(\Omega)}\|\phi_2\|_{H^4(\Omega)} \\
    \leq& C \frac{1+\lambda}{\lambda} \|\phi_1\|_{H^{\max(s,4)}(\Omega)}\|\phi_2\|_{H^{\max(s,4)}(\Omega)}\\
    \leq& C \frac{1+\lambda}{\lambda} R^{\max(s,4)}e^{\sqrt{2}R}. 
    \end{aligned}
\end{equation}
where the first inequality is due to Proposotion~\ref{thm:control}, and the last due to \eqref{eq:CGOHsnorm} .
We obtain the estimate
$$
    \begin{aligned}
    |\hat{\dot{\rho}}(\xi)| \leq&  \left|\int_\Omega\dot{\rho}\phi_1\phi_2\dif x\right| + \left|\int_\Omega\dot{\rho}e^{-i\xi\cdot x}(r_1+r_2+r_1 r_2)\dif x\right|\\
    \leq& C_\lambda\delta + \|\dot{\rho}\|_{H^{-s}(\Omega)}\|r_1+r_2+r_1 r_2\|_{H^s(\Omega)}\\
    \leq& C_\lambda\delta + \|\dot{\rho}\|_{H^{-s}(\Omega)}(\|r_1\|_{H^s(\Omega)}+\|r_2\|_{H^s(\Omega)}+\|r_1\|_{H^s(\Omega)}\|r_2\|_{H^s(\Omega)})\\
    \leq& C_\lambda\delta + C\|\dot{\rho}\|_{H^{-s}(\Omega)}\left(\frac{2}{R}\|q - \lambda\rho_0\|_{H^{s}(\Omega)}+\frac{1}{R^2}\|q - \lambda\rho_0\|_{H^{s}(\Omega)}^2\right)\\
    \leq& C_\lambda\delta + C(\lambda+1)R^{-1}\|\dot{\rho}\|_{H^{-s}(\Omega)},
    \end{aligned}
$$
where in the last inequality we used the estimate~\eqref{eq:rbound}. This derivation holds for any $R\geq 
\frac{|\xi|}{\sqrt{2}}$. In particular, we choose $R=\frac{|\xi|}{\sqrt{2}}$ when $|\xi|>\sqrt{2}a_0(\lambda+1)$, and $R=\sqrt{2}a_0(\lambda+1)$ when $|\xi|\leq\sqrt{2}a_0(\lambda+1)$ to obtain \eqref{eq:fourierest}.

The condition $a_0 \geq C_0 M$ arises since 
$$
C_0\|q-\lambda\rho_0\|_{H^s(\Omega)} \leq C_0 (\|q\|_{H^s(\Omega)} + \lambda \|\rho_0\|_{H^s(\Omega)}) \leq C_0 M (\lambda+1)
$$
is a natural upper bound, thus we require $|\zeta^{(j)}|=R>C_0M(\lambda+1)$ to fulfill the assumption of Lemma~\ref{thm:CGOregularity}.
For either choice of $R$ above, it holds that $R>a_0(\lambda+1)$. It remains to require $a_0\geq C_x0 M$.

\end{proof}

With the help of Lemma~\ref{thm:fourierest}, the following stability estimate can be established for $\dot\rho$.

\medskip   
\begin{thm} \label{thm:stab}
    Let $s>\frac{n}{2}$ with $n\geq 3$. Suppose there exists a constant $M>0$ such that
    $$
    \|q\|_{H^{\max(s,4)}(\Omega)} \leq M, \quad \|\rho_0\|_{H^{\max(s,4)}(\Omega)} \leq M, \quad \|\dot{\rho}\|_{H^s(\Omega)} \leq M.
    $$
    and $\dot\rho$ is compact supported in $\Omega$, then there exist a constant $C$ (independent of $\lambda$ and $\delta$) and a positive constant $\lambda_0>0$ such that
    \[\|\dot\rho\|_{L^\infty(\Omega)}\leq C\left[(\lambda+1)^{\max(s,4)}e^{C(\lambda+1)}\delta+\left(\lambda+\ln\frac{1}{\delta}\right)^{\frac{n-2s}{2}}\right]^\frac{2s-n}{8s}\]
    for any $\lambda>\lambda_0>0$ and $0<\delta\leq e^{-1}$. Here, $e=2.71828...$ is the Euler's number.
\end{thm}

\begin{remark}
For any fixed $\delta>0$, it is clear that $\left(\lambda+\ln\frac{1}{\delta}\right)^{\frac{n-2s}{2}} \rightarrow 0$ as $\lambda\rightarrow\infty$ since $n-2s<0$. Therefore, for a large $\lambda>0$, the estimate in Proposition~\ref{thm:stab} becomes a nearly H\"older-type stability. 
\end{remark}

\begin{proof}
We follow the idea in the proof of the increasing stability result~\cite{nagayasu2013increasing} and name all the constants that are independet of $\lambda$ and $\delta$ as $C$. 

Let $\xi_0$ be a constant such that $\xi_0 \geq \sqrt{2}a_0(\lambda+1)$, then
$$
    \begin{aligned}
    \|\dot{\rho}\|_{H^{-s}(\Omega)}^2=&\int_{\mathbb{R}^n}(1+|\xi|^2)^{-s}|\hat{\dot{\rho}}(\xi)|^2\dif\xi\\
    =&\int_{|\xi|>\xi_0}(1+|\xi|^2)^{-s}|\hat{\dot{\rho}}(\xi)|^2\dif\xi+\int_{\sqrt{2}a_0(\lambda+1)\leq|\xi|\leq \xi_0}(1+|\xi|^2)^{-s}|\hat{\dot{\rho}}(\xi)|^2\dif\xi\\&+\int_{|\xi|\leq \sqrt{2}a_0(\lambda+1)}(1+|\xi|^2)^{-s}|\hat{\dot{\rho}}(\xi)|^2\dif\xi\\
    \eqqcolon& I_1+I_2+I_3.
    \end{aligned}
$$
We estimate $I_1,I_2,I_3$ as follows. For $I_1$, as $\dot\rho$ is compact supported in $\Omega$, H\"older's inequality gives 
$|\hat{\dot\rho}(\xi)|\leq \int_{\Omega}\left|\dot\rho(x)e^{i\xi\cdot x}\right|\dif x\leq C\|\dot\rho\|_{L^2(\Omega)}$. Thus,
\begin{align*}
    I_1 & := \int_{|\xi|>\xi_0}(1+|\xi|^2)^{-s}|\hat{\dot{\rho}}(\xi)|^2\dif\xi \\
    & \leq C\|\dot\rho\|_{L^2(\Omega)}^2\int_{|\xi|>\xi_0}(1+|\xi|^2)^{-s}\dif\xi \\
    & \leq C\|\dot\rho\|_{H^s(\Omega)}^2\xi_0^{n-2s}\leq \underbrace{C \xi_0^{n-2s}}_{:=\Phi_1(\xi_0)} 
\end{align*}
where the last inequality follows from $\|\dot{\rho}\|_{H^s(\Omega)} \leq M$. The function $\Phi_1(\xi_0)$ denotes an upper bound of $I_1$.

For $I_3$, we use $|\hat{\dot{\rho}}(\xi)|\leq \|\hat{\dot{\rho}}\|_{L^\infty(B(0,\sqrt{2}a_0(\lambda+1)))}$ (here $B(0,t)$ means the unit ball of center $0$ and radius $t$) to get
$$
    \begin{aligned}
    I_3 & := \int_{|\xi|\leq \sqrt{2}a_0(\lambda+1)}(1+|\xi|^2)^{-s}|\hat{\dot{\rho}}(\xi)|^2\dif\xi \\
    & \leq \|\hat{\dot{\rho}}\|_{L^\infty(B(0,\sqrt{2}a_0(\lambda+1)))}^2\int_{\mathbb{R}^n}(1+|\xi|^2)^{-s}\dif\xi \\
    & \leq C\|\hat{\dot{\rho}}\|_{L^\infty(B(0,\sqrt{2}a_0(\lambda+1)))}^2  \\
    & \leq C \left[ \frac{(\lambda+1)^2}{\lambda^2}(\lambda+1)^{2\max(s,4)}e^{2\sqrt{2}a_0(\lambda+1)}\delta^2+\frac{1}{a_0^2}\|\dot{\rho}\|_{H^{-s}(\Omega)}^2 \right]. 
    \end{aligned}
$$
where the last inequality is a consequence of~\eqref{eq:fourierest} combined with the estimate $(a+b)^2 \leq 2a^2 + 2b^2$.

For $I_2$, we apply the estimate~\eqref{eq:fourierest} to get
\begin{align*}
I_2 & := \int_{\sqrt{2}a_0(\lambda+1)\leq|\xi|\leq \xi_0}(1+|\xi|^2)^{-s}|\hat{\dot{\rho}}(\xi)|^2\dif\xi \\
& \leq C \int_{\sqrt{2}a_0(\lambda+1)\leq|\xi|\leq \xi_0}(1+|\xi|^2)^{-s}\left|\frac{\lambda+1}{\lambda}|\xi|^{\max(s,4)}e^{|\xi|}\delta+\frac{\lambda+1}{|\xi|}\|\dot{\rho}\|_{H^{-s}(\Omega)}\right|^2\dif\xi\\
& \leq C \int_{\sqrt{2}a_0(\lambda+1)\leq|\xi|\leq \xi_0}(1+|\xi|^2)^{-s}\left[ \frac{(\lambda+1)^2}{\lambda^2}|\xi|^{2\max(s,4)}e^{2|\xi|}\delta^2+\frac{(\lambda+1)^2}{|\xi|^2}\|\dot{\rho}\|^2_{H^{-s}(\Omega)}\right]\dif\xi \\
& = I_{21} + I_{22}
\end{align*}
Let $t:=|\xi|$ be the radial variable, then
\begin{align*}
    I_{21} & := C \int_{\sqrt{2}a_0(\lambda+1)\leq|\xi|\leq \xi_0}(1+|\xi|^2)^{-s} \frac{(\lambda+1)^2}{\lambda^2}|\xi|^{2\max(s,4)}e^{2|\xi|}\delta^2 \dif\xi \\
    & \leq C \frac{(\lambda+1)^2}{\lambda^2} \delta^2 \int^{\xi_0}_{\sqrt{2}a_0(\lambda+1)} (1+t^2)^{-s} t^{2\max(s,4)+n-1} e^{2t} \, dt \\
    & \leq C \frac{(\lambda+1)^2}{\lambda^2} e^{2\xi_0} \delta^2 \int^{\xi_0}_0 t^{2\max(s,4)+n-1-2s} \,dt \\
    & = C \frac{(\lambda+1)^2}{\lambda^2} \xi_0^{2\max(s,4)+n-2s} e^{2\xi_0} \delta^2;
\end{align*}
and
\begin{align*}
    I_{22} & := C \int_{\sqrt{2}a_0(\lambda+1)\leq|\xi|\leq \xi_0}(1+|\xi|^2)^{-s} \frac{(\lambda+1)^2}{|\xi|^2} \|\dot{\rho}\|^2_{H^{-s}(\Omega)} \dif\xi \\
    & = C (\lambda+1)^2 \|\dot{\rho}\|^2_{H^{-s}(\Omega)} \int^{\xi_0}_{\sqrt{2}a_0(\lambda+1)} (1+t^2)^{-s} t^{n-3} \, dt \\
    & \leq C (\lambda+1)^2 \|\dot{\rho}\|^2_{H^{-s}(\Omega)} \int^{\infty}_{\sqrt{2}a_0(\lambda+1)} t^{n-3-2s} \,dt \\
    & \leq C (\lambda+1)^2 \|\dot{\rho}\|^2_{H^{-s}(\Omega)} [\sqrt{2}a_0(\lambda+1)]^{n-2-2s} \\
    & \leq C (\lambda_0+1)^{n-2s} a_0^{n-2s} \frac{1}{a_0^{2}}\|\dot{\rho}\|^2_{H^{-s}(\Omega)} \\
    & = \frac{C}{a_0^{2}}\|\dot{\rho}\|^2_{H^{-s}(\Omega)}.
\end{align*}
Putting these estimates together, we have the following upper bound for $I_2$:
$$
        I_2\leq I_{21} + I_{22}
        \leq \underbrace{ C\frac{(\lambda+1)^2}{\lambda^2}\xi_0^{2\max(s,4)+n-2s}e^{2\xi_0}\delta^2}_{:=\Phi_2(\xi_0)} + \frac{C}{a_0^2}\|\dot{\rho}\|_{H^{-s}(\Omega)}^2.
$$

Combining the estimate for $I_1,I_2,I_3$, we conclude
\begin{equation}
    \begin{aligned}
    \|\dot{\rho}\|_{H^{-s}(\Omega)}^2
    =& I_1 + I_2 + I_3 \\
    \leq & \Phi_1(\xi_0)+ \left[ \Phi_2(\xi_0) + \frac{C}{a_0^2}\|\dot{\rho}\|_{H^{-s}(\Omega)}^2 \right] + \left[C\frac{(\lambda+1)^2}{\lambda^2}(\lambda+1)^{2\max(s,4)}e^{2\sqrt{2}a_0(\lambda+1)}\delta^2+\frac{C}{a_0^2}\|\dot{\rho}\|_{H^{-s}(\Omega)}^2 \right] \\
    \end{aligned}
\end{equation}
where the right hand side has been combined into three groups which are the upper bounds of $I_1,I_2,I_3$, respectively. By choosing $a_0$ sufficiently large, the $H^{-s}$ norm can be absorbed by the left hand side to yield
\begin{equation} \label{eq:rhodotest1}
    \begin{aligned}
    \|\dot{\rho}\|_{H^{-s}(\Omega)}^2\leq& \Phi_1(\xi_0) + \Phi_2(\xi_0) + C\frac{(\lambda+1)^2}{\lambda^2}(\lambda+1)^{2\max(s,4)}e^{2\sqrt{2}a_0(\lambda+1)}\delta^2
    \end{aligned}
\end{equation}

The estimate will henceforth be split into two cases: $\frac{1}{2}\ln\frac{1}{\delta}\geq\sqrt{2}a_0(\lambda+1)$ and $\frac{1}{2}\ln\frac{1}{\delta}<\sqrt{2}a_0(\lambda+1)$. When $\frac{1}{2}\ln\frac{1}{\delta}\geq\sqrt{2}a_0(\lambda+1)$, we choose $\xi_0=\frac{1}{2}\ln\frac{1}{\delta}$ to get
\begin{align*}
\Phi_1(\xi_0) + \Phi_2(\xi_0) = & C\xi_0^{n-2s} + C\frac{(\lambda+1)^2}{\lambda^2}\xi_0^{2\max(s,4)+n-2s}e^{2\xi_0}\delta^2\\
    = & C\left[1+\frac{(\lambda+1)^2}{\lambda^2}\xi_0^{2\max(s,4)}e^{2\xi_0}\delta^2\right]\xi_0^{n-2s} \\
    \leq & C\left[1+\frac{(\lambda+1)^2}{\lambda^2} \left(\ln \frac{1}{\delta} \right)^{2\max(s,4)} \delta \right] \left( \ln\frac{1}{\delta}\right)^{n-2s}.
\end{align*}
As $\lim_{\delta\to0_+}\delta\left(\ln\frac{1}{\delta}\right)^{2\max(s,4)}=0$ and $\lim_{\lambda\rightarrow\infty} \frac{(\lambda+1)^2}{\lambda^2} = 1$, the square parenthesis is bounded whenever $\delta\in(0,e^{-1}]$ and $\lambda\geq \lambda_0$ for some $\lambda_0>0$. 
Hence,
\begin{align*}
\Phi_1(\xi_0) + \Phi_2(\xi_0) & \leq C \left(\ln\frac{1}{\delta}\right)^{n-2s} = C \left(\frac{\ln\frac{1}{\delta}}{\lambda+\ln\frac{1}{\delta}}\right)^{n-2s} \left(\lambda+\ln\frac{1}{\delta}\right)^{n-2s} \\
& \leq C \left(\frac{2\sqrt{2}a_0}{2\sqrt{2}a_0+1}\right)^{n-2s}\left(\lambda+\ln\frac{1}{\delta}\right)^{n-2s} \\
 & \leq C \left(\lambda+\ln\frac{1}{\delta}\right)^{n-2s}
\end{align*}
where the second but last inequality holds since the function $(\frac{t}{\lambda+t})^{n-2s}$ is decreasing in $t>0$.
When $\frac{1}{2}\ln\frac{1}{\delta} < \sqrt{2}a_0(\lambda+1)$, we choose $\xi_0 = \sqrt{2}a_0(\lambda+1)$, then $I_2=0$. As a result, we can simply choose $\Phi_2(\xi_0)=0$ as an upper bound of $I_2$, hence
\begin{align*}
\Phi_1(\xi_0) + \Phi_2(\xi_0) & = \Phi_1(\xi_0) = C \xi_0^{n-2s} = C (\lambda+1)^{n-2s} \\
 & = C \left(\frac{\lambda+1}{\lambda+\ln\frac{1}{\delta}}\right)^{n-2s} \left(\lambda+\ln\frac{1}{\delta}\right)^{n-2s}
 \leq C \left(\frac{1}{1+2\sqrt{2}a_0}\right)^{n-2s}\left(\lambda+\ln\frac{1}{\delta}\right)^{n-2s}.
\end{align*}
In either case, we have
$$
\Phi_1(\xi_0) + \Phi_2(\xi_0) \leq C \left(\lambda+\ln\frac{1}{\delta}\right)^{n-2s}
$$
for some constant $C>0$ that is independent of $\lambda\in [\lambda_0,\infty)$ and $\delta\in (0,e^{-1}]$. In view of~\eqref{eq:rhodotest1}, we conclude
$$
    \begin{aligned}
    \|\dot{\rho}\|_{H^{-s}(\Omega)}^2\leq&C\frac{(\lambda+1)^2}{\lambda^2}(\lambda+1)^{2\max(s,4)}e^{2\sqrt{2}a_0(\lambda+1)}\delta^2+C\left(\lambda+\ln\frac{1}{\delta}\right)^{n-2s}\\
    \leq&C(\lambda+1)^{2\max(s,4)}e^{C(\lambda+1)}\delta^2+C\left(\lambda+\ln\frac{1}{\delta}\right)^{n-2s}.
    \end{aligned}
$$

Finally, we interpolate to obtain an estimate for the infinity norm.
Let $\eta>0$ such that $s=\frac{n}{2}+2\eta$, choose $k_0=-s$, $k_1=s$, $k=\frac{n}{2}+\eta=s-\eta$. Then
\[k=(1-p)k_0+pk_1,\textup{ where }p=\frac{2s-\eta}{2s}.\]
Using the interpolation theorem and the Sobolev embedding, we have
\begin{equation}\label{eq:totalbound}
    \begin{aligned}
    \|\dot{\rho}\|_{L^\infty(\Omega)}\leq& C\|\dot{\rho}\|_{H^k(\Omega)}\leq C\|\dot{\rho}\|_{H^{-s}(\Omega)}^{1-p}\|\dot{\rho}\|_{H^s(\Omega)}^p\leq C\|\dot{\rho}\|_{H^{-s}(\Omega)}^{\frac{2s-n}{8s}}\\
    \leq& C\left[(\lambda+1)^{2\max(s,4)}e^{C(\lambda+1)}\delta^2+C\left(\lambda+\ln\frac{1}{\delta}\right)^{n-2s}\right]^\frac{2s-n}{16s}\\
    \leq& C\left[(\lambda+1)^{\max(s,4)}e^{C(\lambda+1)}\delta+C\left(\lambda+\ln\frac{1}{\delta}\right)^{\frac{n-2s}{2}}\right]^\frac{2s-n}{8s}\\
    \end{aligned}
\end{equation}

\end{proof}

\section{Numerical Experiments} \label{sec:experiment}

\subsection{Setup}
In this section we implement Algorithm~\ref{alg:constq} for the case $\rho_0=1$, $q\equiv 0$ in one dimension(1D). 
We choose $\Omega=[-1,1]$ as the computational domain, $T=5$ as the terminal time, and $\lambda = \lambda_j = \frac{j^2 \pi^2}{4}$ with $j=1,2,\dots,10$. The idea is to compute the Fourier coefficients of $\dot\rho$ using~\eqref{eq:keyid} with respect to the Fourier basis functions
\begin{equation} \label{eq:Fourierbasis}
    \left\{1, \; \cos(j\pi x), \; \sin(j\pi x) \right\}^{10}_{j=1}.
\end{equation}

Specifically, for each $\lambda_j>0$, we solve the boundary control equations
$$
u^f_0(T) = \cos(\frac{j\pi}{2}x), \quad\quad
u^h_0(T) = \cos(\frac{j\pi}{2}x),
$$
then the Fourier coefficient $(\dot{\rho},\cos(j\pi x))_{L^2(-1,1)}$ can be computed using~\eqref{eq:keyid}. If we instead solve the boundary control equations
$$
u^f_0(T) = \cos(\frac{j\pi}{2}x), \quad\quad
u^h_0(T) = \sin(\frac{j\pi}{2}x),
$$
then the Fourier coefficient $(\dot{\rho},\sin(j\pi x))_{L^2(-1,1)}$ becomes available. 
It remains to obtain the zeroth-order Fourier coefficient $(\dot{\rho},1)_{L^2(-1,1)}$. This term, however, cannot be computed by solving the boundary control equations
$u^f_0(T) = u^h_0(T) = 1$, since these Helmholtz solutions correspond to the eigenvalue $\lambda=0$, for which the left hand side of~\eqref{eq:keyid} vanishes. Instead, we can take an arbitrary positive eigenvalue $\lambda_j = \frac{j^2\pi^2}{4}$ for some $j$, compute the inner products $(\dot\rho, \cos^2(\frac{j\pi}{2}x))_{L^2(-1,1)}$ and $(\dot\rho, \sin^2(\frac{j\pi}{2}x))_{L^2(-1,1)}$ using~\eqref{eq:keyid}, then add them to get $(\dot{\rho},1)_{L^2(-1,1)}$.

\subsection{Solving the Boundary Control Equations} \label{sec:TR}
Here, we adopt the analytic approach in the authors' earlier work~\cite{oksanen2022linearized} to solve the boundary control equations. The idea is to back-propagate a desired final state $\phi(x)$ to $t=0$ in $\mathbb{R}^n$, then the Neumann trace is a suitable boundary control.

Specifically, for $\phi\in C^\infty([-1,1])$, define \[
\tilde\phi := \begin{cases}
    \phi&x\in[-1,1],\\
    \phi\cdot\exp\{1-\frac{1}{1-(x+1)^{4}}\}&x\in(-2,-1),\\
    \phi\cdot\exp\{1-\frac{1}{1-(x-1)^{4}}\}&x\in(1,2),\\
    0&x\notin(-2,2),
    \end{cases}
\]
then $\tilde{\phi}$  is a $C^3$-extension of $\phi$. According to the D'Alembert formula, the solution of the problem
\[
\left\{
\begin{alignedat}{2}
    u_{tt}-u_{xx} &= 0\\
    u(0,x) = u_t(0,x) &= 0\\
    u(T,x) &= \phi(x)
\end{alignedat}
\right.
\]
is 
\[u(t,x) = \frac{1}{2}[\tilde\phi(x+t-T)+\tilde\phi(x-t+T)]\]
when $T$ is sufficiently large. Thus, if we choose
\[f(t,\pm 1) := \left.\pm\frac{1}{2}[\tilde\phi'(x+t-T)+\tilde\phi'(x-t+T)]\right|_{x=\pm1},\]
then $u_0^f(T) = \phi$. In this case, $\partial^2_t f$ can be analytically computed.

\subsection{Discretization}
Since the right hand side of~\eqref{eq:linearizedBVP} involves the second order derivative $\partial^2_t u_0$, we use a fourth order finite difference method for the spatial discretization so that the discretization error of the second derivative converges to zero as the grid size tends to zero.
For a function $f$, the finite difference approximation of its $k$-th order derivative using $m$ ($k \leq m-1$) grid points is 
\[f^{(k)}(x)\approx \sum_{i=1}^mc_if(x+h_i).\]
Matching the Taylor expansions up to order $m-1$ shows that the coefficients $c_i$ should satisfy 
\begin{equation}\label{eq:fdmcoef}
\mathbf{A}\mathbf{c}=k!\times \mathbf{e}_k
\end{equation}
to guarantee at least $(m-k)$-th order accuracy. Here, $\mathbf{A}$ is an $m\times m$ Vandermonde matrix with entries $\mathbf{A}_{ij}=h_j^{i-1}$; $\mathbf{c}$ is the vector whose components are the coefficients $c_i$; $\mathbf{e}_k$ is the vector whose $k$-th component is one and all the others are zeros.

We choose the spatial grid points $x_j:= -1+j \Delta x$ with spacing $\Delta x$, $j=0,1,\dots,N$. In the following experiments, we choose $\Delta x = 0.004$ so that $N=500$. To apply the Neumann derivative, we add two additional ghost points $x_{-1}:=x_0-\Delta x$, $x_{N+1}:=x_N+\Delta x$. We choose $T=5$ with $\Delta t=0.1\Delta x$ to fulfill the CFL condition.
Solving \eqref{eq:fdmcoef} gives the following fourth order finite difference approximations of the first order derivatives as well as the second order derivatives:
\begin{equation}\label{eq:4thfdm}
    \begin{aligned}
        f'(x_0)&\approx\frac{1}{12\Delta x}(-3f(x_{-1})-10f(x_0)+18f(x_1)-6f(x_2)+f(x_3)),\\
        f'(x_N)&\approx\frac{1}{12\Delta x}(3f(x_{N+1})+10f(x_N)-18f(x_{N-1})+6f(x_{N-2})-f(x_{N-3})),\\
        f''(x_0)&\approx\frac{1}{12\Delta x^2}(10f(x_{-1})-15f(x_0)-4f(x_1)+14f(x_2)-6f(x_3)+f(x_4)),\\
        f''(x_i)&\approx\frac{1}{12\Delta x^2}(-f(x_{i-2})+16f(x_{i-1})-30f(x_i)+16f(x_{i+1})-f(x_{i+2})),\qquad(1\leq i \leq N-1)\\
        f''(x_N)&\approx\frac{1}{12\Delta x^2}(f(x_{N-4})-6f(x_{N-3})+14f(x_{N-2})-4f(x_{N-1})-15f(x_N)+10f(x_{N+1})),\\
    \end{aligned}    
\end{equation}

The wave equation is solved using the finite difference time domain method. For each time step, the values of the solution at $x_{-1}$ and $x_{N+1}$ are updated using the first two approximations in \eqref{eq:4thfdm} along with the boundary condition $f'(x_0)=f'(x_N)=0$, then we calculate $\Delta u$ and $u_{tt}$ using the last three approximations in \eqref{eq:4thfdm}. The next time step is updated using $u_{tt}$ and the last approximation in \eqref{eq:4thfdm}. Since the Neumann data for each wave equation close to $t=0$ is zero and the initial conditions are zeros, the numerical solution for the first few time steps are zero when $\Delta t$ is sufficiently small, we can start iteration at the sixth time step to avoid involving wave solution before initial state.

\subsection{Numerical Experiment}

Recall that $\rho_0\equiv 1$ and $q\equiv 0$ in all the experiments.

\begin{figure}[!h]
    \centering
    \includegraphics[width=0.49\textwidth]{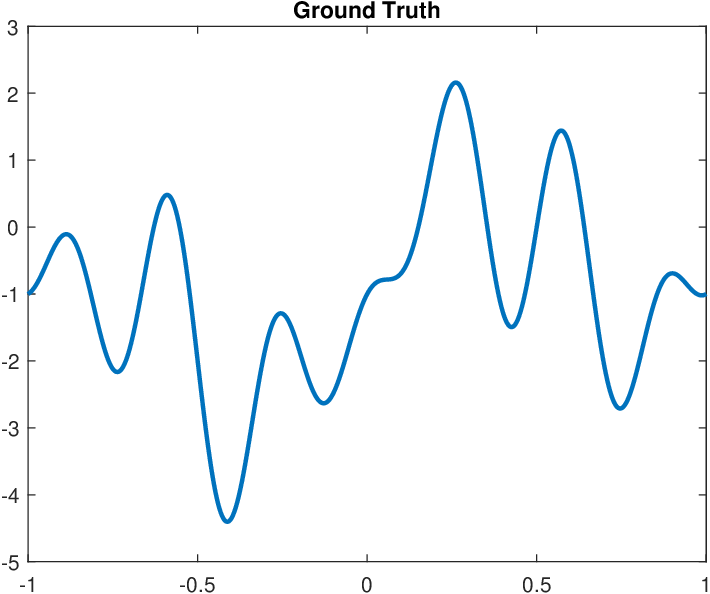}
    \caption{Ground truth $\dot\rho$}
    \label{fig:groundtruth1}
\end{figure}

\begin{figure}[!h]
    \centering
    \includegraphics[width=0.49\textwidth]{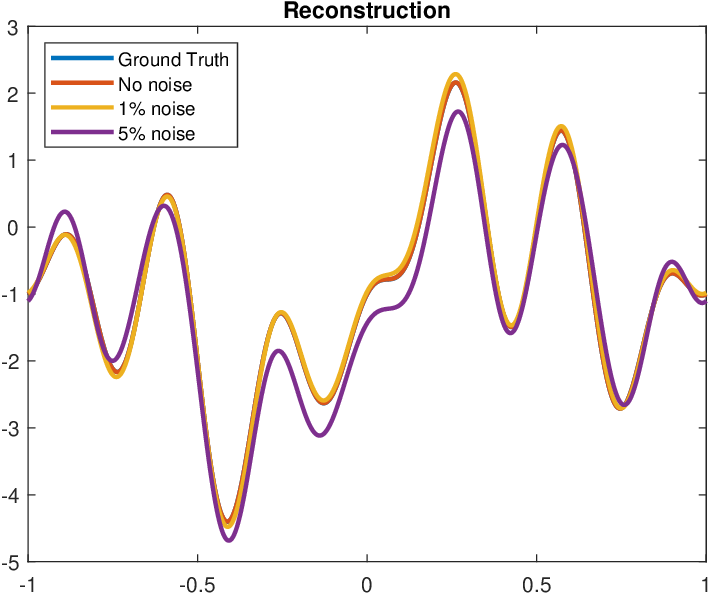}
    \includegraphics[width=0.49\textwidth]{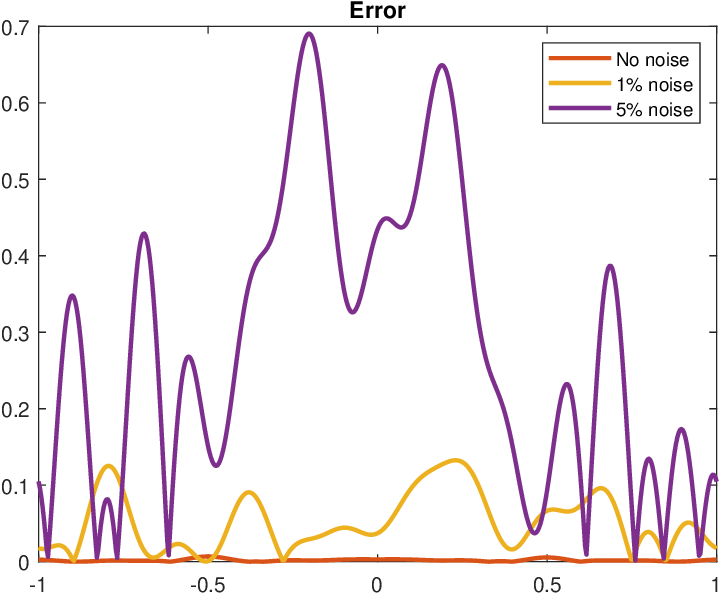}
    \caption{Left: Reconstructed $\dot\rho$ with $0\%,1\%,5\%$ Gaussian noise and the ground truth. Right: The corresponding error between the reconstruction result and the ground truth. The relative $L^2$-errors are $0.14\%,3.66\%$ and $19.37\%$, respectively.}
    \label{fig:reconstruction1}
\end{figure}

\textbf{Experiment 1.}
We start with a continuous perturbation 
\[\dot{\rho} = \sin(\pi x) + \sin(2\pi x) - \cos(5\pi x) + \cos(7\pi x) - 1,\]
which is in the span of the Fourier basis functions~\eqref{eq:Fourierbasis}. The graph of $\dot\rho$ is shown in Figure~\ref{fig:groundtruth1}. The Gaussian random noise are added to the measurement $\dot\Lambda$ by adding to the numerical solutions on the boundary nodes. The reconstructions and corresponding errors 
 with noise level $0\%,1\%,5\%$ are illustrated in Figure~\ref{fig:reconstruction1}.

\begin{figure}[!h]
    \centering
    \includegraphics[width=0.49\textwidth]{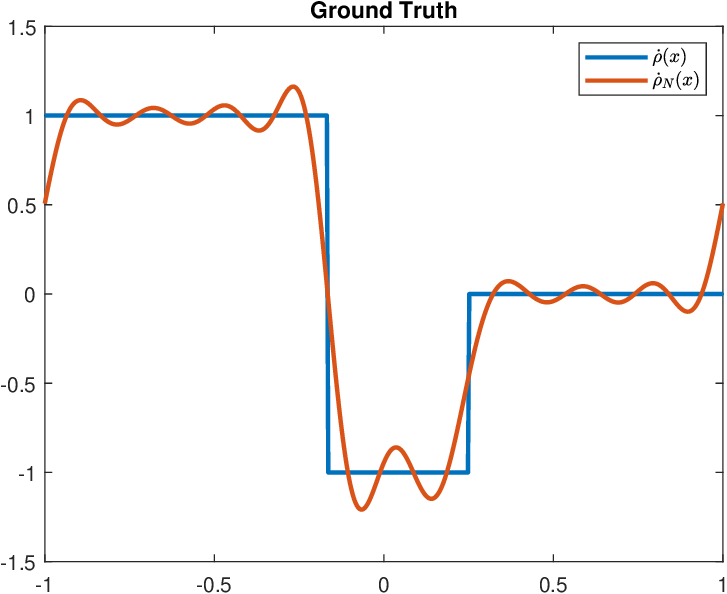}
    \caption{Ground truth $\dot\rho$}
    \label{fig:groundtruth2}
\end{figure}

\begin{figure}[!h]
    \centering
    \includegraphics[width=0.49\textwidth]{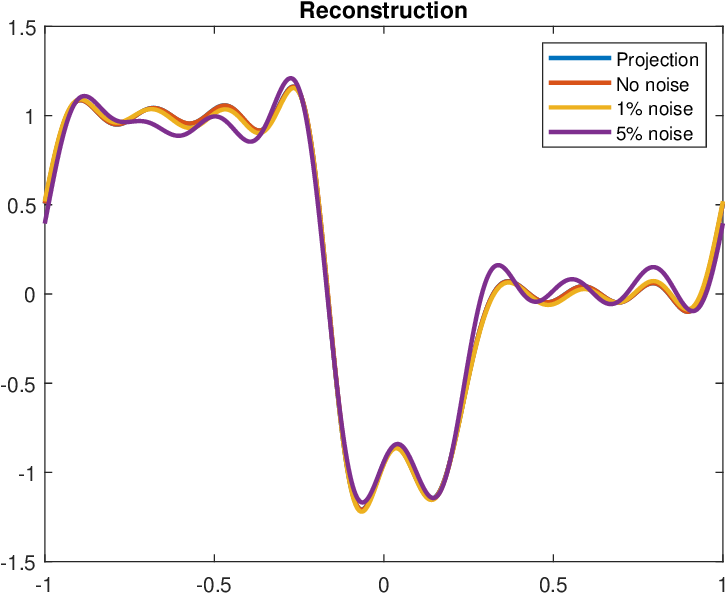}
    \includegraphics[width=0.49\textwidth]{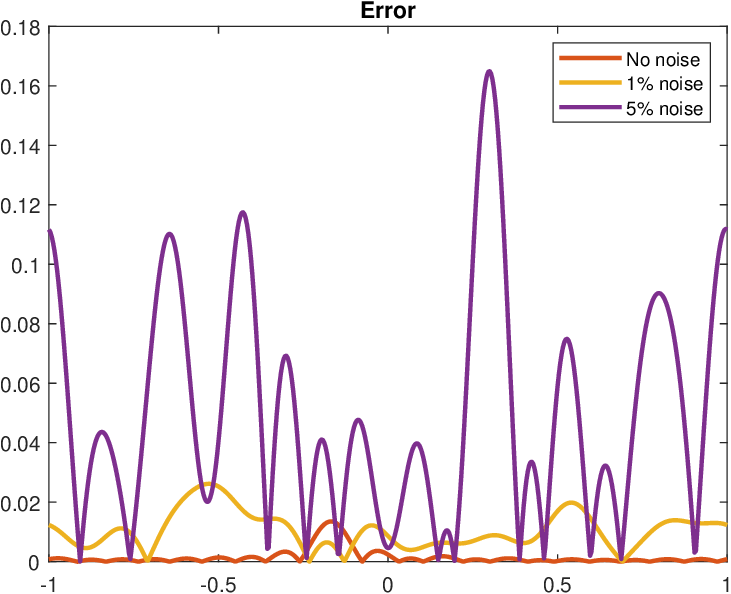}
    \caption{Left: Reconstructed $\dot\rho$ with $0\%,1\%,5\%$ Gaussian noise and the ground truth. Right: The corresponding error between the reconstruction result and the ground truth. The relative $L^2$-errors are $0.39\%,1.57\%$ and $8.15\%$, respectively.}
    \label{fig:reconstruction2}
\end{figure}

\textbf{Experiment 2.}
In this experiment, we consider a discontinuous perturbation
\[\dot\rho=\chi_{[-1,-\frac{1}{6}]}-\chi_{[-\frac{1}{6},\frac{1}{4}]},\]
where $\chi$ is the characteristic function. The Fourier series of $\dot\rho$ is given by
\[\dot{\rho}=\frac{5}{24}+\sum_{n=1}^\infty\left[-\frac{\sin\left(\frac{n\pi}{4}\right)+2\sin\left(\frac{n\pi}{6}\right)}{n\pi}\cos(n\pi x)+\frac{\cos(n\pi)+\cos\left(\frac{n\pi}{4}\right)+2\cos\left(\frac{n\pi}{6}\right)}{n\pi}\sin(n\pi x)\right].\]
With the choice of the basis functions~\eqref{eq:Fourierbasis}, we can only expect to reconstruct the orthogonal projection:
\[\dot{\rho}_N\coloneqq\frac{5}{24}+\sum_{n=1}^N\left[-\frac{\sin\left(\frac{n\pi}{4}\right)+2\sin\left(\frac{n\pi}{6}\right)}{n\pi}\cos(n\pi x)+\frac{\cos(n\pi)+\cos\left(\frac{n\pi}{4}\right)+2\cos\left(\frac{n\pi}{6}\right)}{n\pi}\sin(n\pi x)\right],\]
see Figure~\ref{fig:groundtruth2}. We plot the reconstruction result and the corresponding error with respect to the orthogonal projection $\dot{\rho}_N$ with different noise level in Figure~\ref{fig:reconstruction2}.


\begin{figure}[!h]
    \centering
    \includegraphics[width=0.49\textwidth]{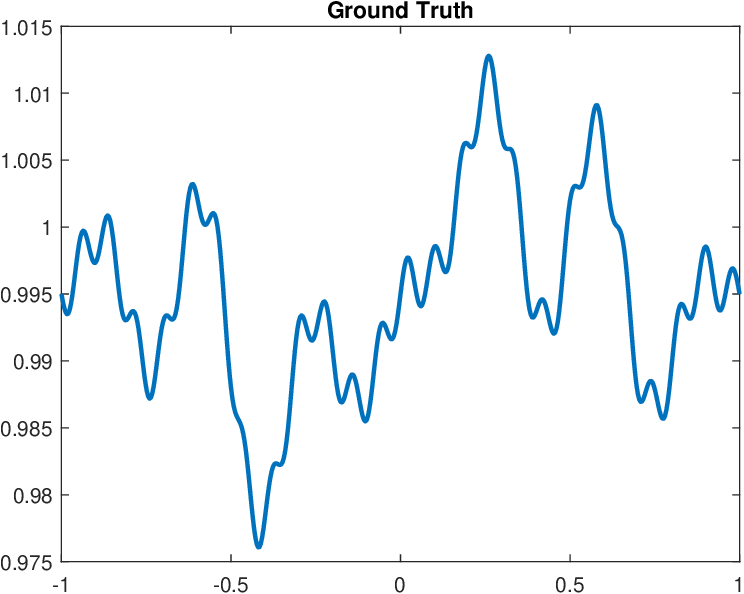}
    \caption{Ground truth $\rho$}
    \label{fig:groundtruth3}
\end{figure}

\begin{figure}[!h]
    \centering
    \includegraphics[height=0.25\textheight]{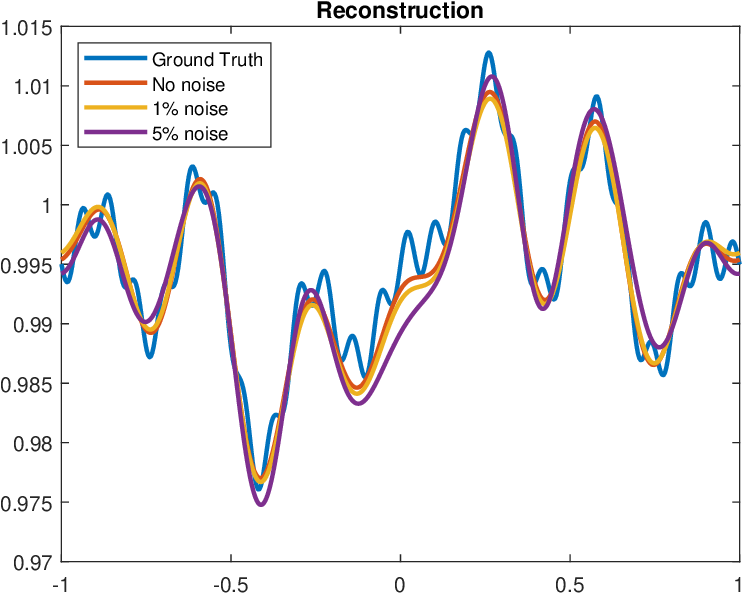}
    \includegraphics[height=0.25\textheight]{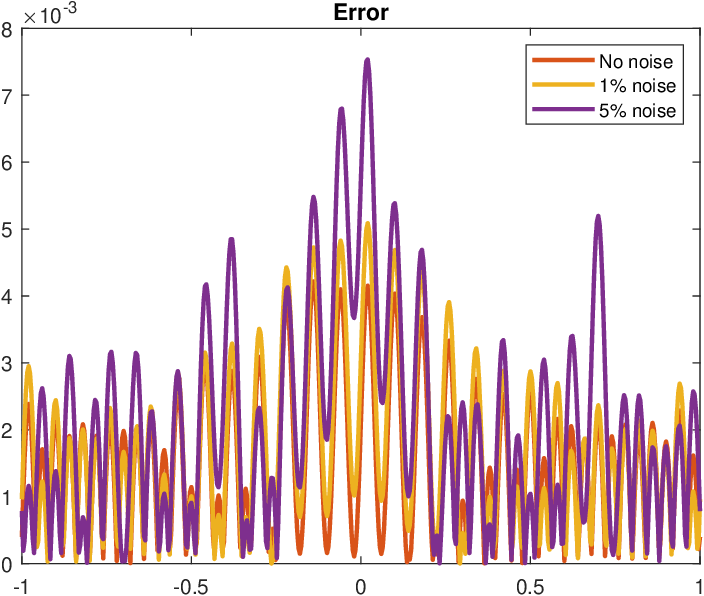}
    \caption{Left: Reconstructed $\rho$ with $0\%,1\%,5\%$ Gaussian noise and the ground truth. Right: The corresponding error between the reconstruction result and the ground truth. The relative $L^2$-errors are $18.09\%,20.95\%$ and $26.46\%$, respectively.}
    \label{fig:reconstruction3}
\end{figure}

\textbf{Experiment 3.}
In this experiment, we apply the algorithm to the non-linear IBVP where
\[\rho = \rho_0 + \varepsilon\dot{\rho} + \varepsilon^2\ddot{\rho},\]
with $\varepsilon=0.001$ and
\[\dot{\rho} = \sin(\pi x) + \sin(2\pi x) - \cos(5\pi x) + \cos(7\pi x) - 1,\qquad\ddot{\rho} = 200\sin(25\pi x).\]
See Figure~\ref{fig:groundtruth3} for the graph of $\rho$.
Since
\[\Lambda_\rho-\Lambda_{\rho_0}\approx\varepsilon\dot\Lambda_{\dot\rho}=\dot\Lambda_{\varepsilon\dot\rho}\]
when $\varepsilon$ is small, we can use $\varepsilon^{-1}(\Lambda_\rho-\Lambda_{\rho_0})$ as an approximation of $\dot\Lambda_{\dot\rho}$ in \eqref{eq:keyid}. In this case, $\Lambda_\rho f$ and $\Lambda_{\rho_0}f$ are computed by numerically solving the forward problem \eqref{eq:unperturbedBVP} with $\rho$ and $\rho_0\equiv 1$. We then apply Algorithm~\ref{alg:constq} to find $\dot{\rho}$, and view $1+\dot\rho$ as an approximation of the ground truth $\rho$.
In the experiment, we added the Gaussian noise to the difference $\Lambda_\rho-\Lambda_{\rho_0}$ rather than to $\Lambda_\rho f$ and $\Lambda_{\rho_0}f$ individually, see~\cite{oksanen2022linearized} for discussion of the difference. The reconstruction and the respective errors with $0\%,1\%,5\%$ Gaussian noise are illustrated in Figure~\ref{fig:reconstruction3}.


\appendix

\subsection*{Appendix}

We provide a rigorous proof to show that the operator $\dot{\Lambda}_{\dot{\rho}}$ derived in the introduction using the formal argument is indeed the Fr\'echet derivative of the nonlinear map $\rho\mapsto \Lambda_\rho$ at $\rho_0$ along $\dot{\rho}$.

Denote by $\mathcal{L}(X,Y)$ the Banach space of bounded linear operators from $X$ to $Y$. The IBVP~\eqref{eq:bvp} aims to invert the following nonlinear map
\[\mathcal{F}:\rho\in C^\infty(\overline{\Omega})\mapsto\Lambda_\rho\in\mathcal{L}(H^{\frac{5}{2}}((0,2T)\times\partial\Omega),H^{\frac{1}{2}}((0,2T)\times\partial\Omega))\]
Suppose $\rho = \rho_0 + \dot\rho$ with $\rho_0\in C^\infty(\overline{\Omega})$ and $\dot{\rho}\in C_c^\infty(\Omega)$, define the following linear operator
\[\dif\mathcal{F}:\dot\rho\in C_c^\infty(\Omega)\mapsto\dot{\Lambda}_{\dot{\rho}}\in\mathcal{L}(H^{\frac{5}{2}}((0,2T)\times\partial\Omega),H^{\frac{1}{2}}((0,2T)\times\partial\Omega))\]
where $\dot{\Lambda}_{\dot{\rho}}$ is the linearized ND map defined in \eqref{eq:linearizedNDmap}.

\begin{prop}
    The nonlinear map $\mathcal{F}$ is Frech\'et differentiable at $\rho_0\in C^\infty(\overline{\Omega})$, and the Frech\'et derivative along the direction $\dot\rho\in C_c^\infty(\Omega)$ is $\dot{\Lambda}_{\dot{\rho}}$.
\end{prop}
\begin{proof}
In order to show that $\mathcal{F}$ is Frech\'et differentiable, we will show that
\[\|\mathcal{F}(\rho)-\mathcal{F}(\rho_0)-\dif\mathcal{F}(\dot\rho)\|_{\mathcal{L}(H^{\frac{5}{2}}((0,2T)\times\partial\Omega),H^{\frac{1}{2}}((0,2T)\times\partial\Omega))} = O(\|\dot\rho\|_{W^{1,\infty}(\Omega)}^2)\]
as $\|\dot\rho\|_{W^{1,\infty}(\Omega)}\to 0$, which is equivalent to
\[\|\Lambda_\rho f - \Lambda_{\rho_0} f - \dot{\Lambda}_{\dot\rho} f\|_{H^{\frac{1}{2}}((0,2T)\times\partial\Omega)} =O( \|\dot\rho\|_{W^{1,\infty}(\Omega)}^2\|f\|_{H^{\frac{5}{2}}((0,2T)\times\partial\Omega)})\]
for any $f\in H^{\frac{5}{2}}((0,2T)\times\partial\Omega)$ as $\|\dot\rho\|_{W^{1,\infty}(\Omega)}\to 0$.

Write $u=u_0+\delta u$, where $u$ and $u_0$ are the solutions of \eqref{eq:bvp} and \eqref{eq:unperturbedBVP}, respectively. Then $\delta u$ satisfies the equation 
$$
\left\{
\begin{alignedat}{2}
    \rho_0\delta u_{tt} - \Delta\delta u + q\delta u &= -\dot{\rho} \partial^2_t u\qquad\qquad&&\text{in }(0,2T)\times\Omega\\
    \partial_\nu \delta u &= 0&&\text{on }(0,2T)\times\partial\Omega\\
    \delta u(0,x) = \delta u_t(0,x) &= 0&& x\in\Omega
\end{alignedat}
\right.
$$
Using the regularity estimate for wave euation and the trace theorem, we have
$$
    \|\delta u\|_{H^2((0,2T)\times\Omega)}\leq C\|\dot{\rho} u_{tt}\|_{H^1((0,2T)\times\Omega)}\leq C\|u\|_{H^3((0,2T)\times\Omega)}\|\dot\rho\|_{W^{1,\infty}(\Omega)}
$$

On the other hand, denote $w\coloneqq \delta u - \dot{u}$; it satisfies 
$$
\left\{
\begin{alignedat}{2}
    \rho_0w_{tt} - \Delta w + qw &= -\dot{\rho} \partial^2_t \delta u \qquad\qquad&&\text{in }(0,2T)\times\Omega\\
    \partial_\nu w &= 0&&\text{on }(0,2T)\times\partial\Omega\\
    w(0,x) = w_t(0,x) &= 0&& x\in\Omega
\end{alignedat}
\right.
$$
Applying similar estimates yield
$$
    \|w\|_{H^\frac{1}{2}((0,2T)\times\partial\Omega)}\leq\|w\|_{H^1((0,2T)\times\Omega)}\leq C\|\dot{\rho} \delta u_{tt}\|_{L^2((0,2T)\times\Omega)}\leq C\|\delta u\|_{H^2((0,2T)\times\Omega)}\|\dot\rho\|_{W^{1,\infty}(\Omega)}.
$$
Notice that $w|_{[0,2T]\times\partial\Omega} = \Lambda_\rho f - \Lambda_{\rho_0} f - \dot{\Lambda}_{\dot\rho} f$. These estimates combined give
$$
\|\Lambda_\rho f - \Lambda_{\rho_0} f - \dot{\Lambda}_{\dot\rho} f \|_{H^\frac{1}{2}((0,2T)\times\partial\Omega)} \leq C\|u\|_{H^3((0,2T)\times\Omega)}\|\dot\rho\|^2_{W^{1,\infty}(\Omega)}.
$$
Finally, the well-posedness estimate for the forward boundary value problem~\cite{evans1997pde} gives 
$$
\|u\|_{H^3((0,2T)\times\Omega)} \leq C \|f\|_{H^{\frac{5}{2}}((0,2T)\times\partial\Omega)}.
$$
This completes the proof.
\end{proof}


\section*{Acknowledgement}
The research of T. Yang and Y. Yang is partially supported by the NSF grants DMS-2006881, DMS-2237534, DMS-2220373, and the NIH grant R03-EB033521. L. Oksanen is supported by the European Research Council of the European Union, grant 101086697 (LoCal), and the Reseach Council of Finland, grants 347715 and 353096. Views and opinions expressed are those of the authors only and do not necessarily reflect those of the European Union or the other funding organizations.


\bibliographystyle{abbrv}
\bibliography{refs}

\end{document}